\documentclass[12pt]{article}
\usepackage{amssymb}
\usepackage{amsthm}
\theoremstyle{plain}
\newtheorem{tetel}{Theorem}[section]
\newtheorem{all}[tetel]{Proposition}
\newtheorem{lemma}[tetel]{Lemma}
\newtheorem{kov}[tetel]{Corollary}

\theoremstyle{definition}\newtheorem{Def}[tetel]{Definition}

\theoremstyle{remark}\newtheorem{megj}[tetel]{Remark}
\newtheorem{pelda}[tetel]{Example}

\usepackage{amsmath}
\usepackage[latin2]{inputenc}
\usepackage{t1enc}
\usepackage[mathscr]{eucal}
\usepackage{indentfirst}
\usepackage{graphicx}
\usepackage{psfrag}
\usepackage{fancyhdr}
\usepackage{amscd}
\usepackage{enumerate}
\usepackage{floatflt}
\newcommand*{\R}{\ensuremath{\mathbf R}}

\newcommand*{\Z}{\ensuremath{\mathbf Z}}

\DeclareMathOperator{\id}{id}

\DeclareMathOperator{\lcm}{lcm}

\usepackage{array}
\usepackage{tabularx}

\begin{document}

\title{One parameter families of Legendrian torus knots}
\author{Tam\'as K\'alm\'an\\
\small University of Southern California}
\date{}
\maketitle

\begin{abstract}
We compute the Chekanov--Eliashberg contact 
homology of what we call the Legendrian closure of a positive braid 
(Definition \ref{def:lezaras}). 
We also construct an augmentation for each such link diagram. 
Then we apply the techniques established in the first of this series of papers 
\cite{en} to a certain natural loop in the 
space $\mathscr L'$ of positive Legendrian $(p,q)$ torus knots to 
show that $\mathscr L'$ has a nontrivial 
fundamental group, which is mapped non-injectively into $\pi_1(\mathscr 
K)$ by the map induced by inclusion (here, 
$\mathscr K$ is the space of all smooth positive $(p,q)$ torus knots). 
This is the first example that Legendrian and classical knots 
behave differently at the level of one parameter families. 
\end{abstract}

\section{Introduction}

We assume that the reader is familiar with 
\cite{en}. Note however that 
in this paper, we work in the original, $\Z_2$--coefficient 
version of Chekanov--Eliashberg theory; thus the formulas of 
Table 1, Theorem 3.2, and Remark 3.4 
of \cite{en} have to be simplified accordingly. Namely, in the knot DGA, 
we set $t=1$ and 
reduce all coefficients modulo $2$, and we reduce the index of each 
generator modulo $2r$. Thus the grading of the contact homology also 
becomes a $\Z_{2r}$--grading. In section \ref{sec:z2}, we review the 
chain maps that define monodromies of loops. All the relevant theorems of 
\cite{en} hold in this reduced setting. 

The main purpose of this paper is to prove Theorem \ref{thm:monorendje} 
(which almost immediately implies Theorem \ref{thm:hurokrendje} and 
Corollary \ref{kov}) about torus knots. Below we give an outline of the 
proof. However in sections \ref{sec:BC}, \ref{sec:aug}, and 
\ref{sec:loop}, we will work in the more general case of positive 
links, which turn out to have a natural Legendrian representative (see 
Figure \ref{fig:pic}) with remarkable properties (see especially 
Remark \ref{rem:fuchs}). 
While we hope that the concluding sections \ref{sec:moreBC} and 
\ref{sec:moreaug} 
(i.e., the computation of the order of the monodromy, see below) 
will also be suitably generalized to arbitrary positive links in future work, 
we do believe that Theorem \ref{thm:relacio} (contact homology) and 
Proposition 
\ref{pro:aug} (augmentation) are interesting on their own right, too.

Let us consider the space $\mathscr L'$ of positive Legendrian 
$(p,q)$ torus knots, where $p,q\ge2$ are relatively prime integers. 
In \cite{EH1}, Etnyre and Honda showed that the path-components  
of $\mathscr L'$ are completely classified by their Thurston--Bennequin 
number $tb$ and rotation number $r$ (i.e.\ that $\mathscr K$, the space 
of 
smooth positive $(p,q)$ torus knots, is \emph{Legendrian simple}). 
In particular, according to their classification, there is a single 
component $\mathscr L$ to which the maximal value $tb=pq-p-q$ is 
associated (this one has rotation number $r=0$), 
and all other positive Legendrian 
$(p,q)$ torus knots are stabilizations\footnote{A stabilization is the result 
of performing a Reidemeister I move on a Lagrangian diagram.} of this type.

In section \ref{sec:BC}, we compute the Chekanov--Eliashberg 
contact homology $H$ 
of $\mathscr L$, i.e.\ we find all admissible discs in a Lagrangian diagram 
$\gamma$ representing a certain knot $L\in\mathscr L$. In section 
\ref{sec:aug}, we introduce an augmentation $\varepsilon$ of $\gamma$. (An 
augmentation can be thought of as a $\Z_2$--valued cocycle on the DGA which 
is also a graded algebra morphism.) In 
section \ref{sec:loop}, we define an 
$S^1$--family $\Omega_{p,q}\subset\mathscr L$ based at $L$. Let us 
denote its monodromy by $\mu\colon H(L)\to H(L)$ and let 
$\mu_0\colon H_0(L)\to H_0(L)$ be the 
restriction of $\mu$ to the index $0$ part of $H(L)$. 
In section \ref{sec:moreBC}, 
we use a direct computation to show that 
$\mu_0^{p+q}=\id_{H_0(L)}$, i.e.\ we show that the order of $\mu_0$ divides 
$p+q$. 
This by itself implies nothing for topology, but the formulas we develop are 
needed in section \ref{sec:moreaug}, where we use $\varepsilon$ to prove 
that the order of $\mu_0$, and thus that of $\mu$, 
is divisible by $p+q$. By accomplishing all this, 
we will have established the following:

\begin{tetel}\label{thm:monorendje}
The restricted monodromy $\mu_0$ of the loop $\Omega_{p,q}$ has order $p+q$.
\end{tetel}

In particular, $\mu$ is different from the identity automorphism. Note again 
that the proof of this claim depends heavily on the augmentation that we 
construct in section \ref{sec:aug}.

Another observation about $\Omega_{p,q}$, 
based on results of Goldsmith \cite{gold} 
(which we also noted in the introduction of 
\cite{en}), is that it has order $2p$ in $\pi_1(\mathscr K,L)$. These 
facts combined with Theorem 1.1 
of \cite{en} yield that

\begin{tetel}\label{thm:hurokrendje}
The order of $[\Omega_{p,q}]\in\pi_1(\mathscr L,L)$ 
is either infinite or divisible 
by $\lcm\{\,2p,p+q\,\}$ (which is either $p(p+q)$ or $2p(p+q)$).
\end{tetel}

Let us point out again that by the theorem, 
$[\Omega_{p,q}]^{2p}\in\pi_1(\mathscr L,L)$ 
is nontrivial, which implies the 

\begin{kov}\label{kov}
There exist Legendrian knot types $\mathscr L$ so that for the
corresponding smooth knot type $\mathscr K\supset\mathscr L$, the
homomorphism $\pi_1(\mathscr L)\to\pi_1(\mathscr K)$ induced by the
inclusion is not injective.
\end{kov}

\bigskip

Acknowledgements: This is the second half of the author's thesis that was 
submitted in 2004 at the University of California, Berkeley. 
I would like to thank Michael Hutchings, Rob Kirby, 
Lenny Ng, Josh Sabloff, Michael   
Sullivan and Andr\'as Sz\H{u}cs for stimulating conversations. I was
supported by NSF grant DMS-0244558 during part of the research.

\section{\boldmath Monodromy maps with 
$\Z_2$--coefficients}\label{sec:z2}

\begin{figure}[t]
\psfrag{a}{\Huge $a$}
\psfrag{b}{\Huge $b$}
\psfrag{c}{\Huge $c$}
\psfrag{3a}{\Huge{III$_{\text{a}}$}}
\psfrag{3b}{\Huge{III$_{\text{b}}$}}
\psfrag{-1}{\Huge{II}}
\psfrag{2}{\Huge{II$^{-1}$}}
\psfrag{a'}{\Huge $a'$}
\psfrag{b'}{\Huge $b'$}
\psfrag{c'}{\Huge $c'$}
\resizebox{\linewidth}{!}{\includegraphics{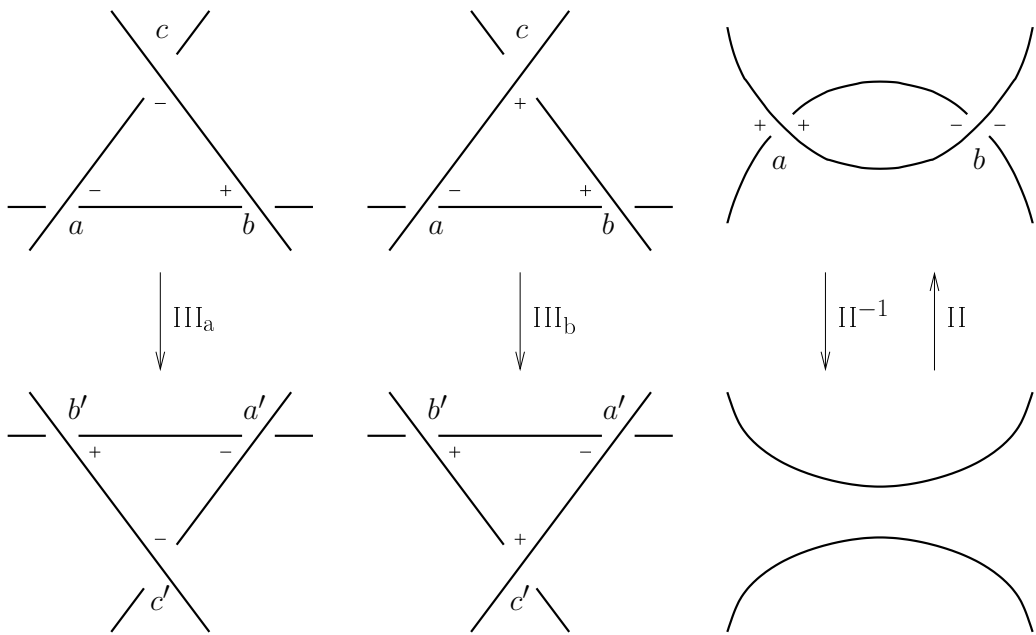}}
\caption{Legendrian Reidemeister moves (of Lagrangian projections). The 
signs shown are Reeb signs.}\label{fig:moves}
\end{figure}

This is a summary of section 3.1 
of \cite{en}, reduced to the original $\Z_2$--coefficient theory of 
Chekanov. I.e., we give a list of maps that, when extended from the 
generators to the DGA as algebra 
morphisms, become the chain maps that are used to define 
holonomies and monodromies of (sequences of) Reidemeister moves. 

Please refer to Figure \ref{fig:moves} and recall that in \cite{en} we 
denoted the DGA's of the
diagrams in the upper row by $(\mathscr A,\partial)$ and in the lower
row, by $(\mathscr A',\partial')$. The crossings not affected by the
moves are in an obvious one-to-one correspondence. If $x\in\mathscr A$
is one of these, then the corresponding crossing in $\mathscr A'$ 
is denoted by $x'$.

Move III$_{\text{a}}$: 
Let $a\mapsto a'$, $b\mapsto b'$, $c\mapsto c'$, and $x\mapsto x'$,
where $x$ is any other crossing of the upper diagram.

Move III$_{\text{b}}$: 
Let 
\[a\mapsto a'+c'b',\] 
while other generators are
mapped trivially: $b\mapsto b'$, $c\mapsto c'$, and $x\mapsto x'$.

Move II$^{-1}$:
Let $\partial(a)=b+v$ (recall that $v$ does not contain neither
$a$ nor $b$). Define $x\mapsto x'$, which gives rise to the obvious
re-labeling $v\mapsto v'$. Then let
\begin{eqnarray*}
a&\mapsto&0\\
b&\mapsto&v'.
\end{eqnarray*}

Move II: 
Suppose that right after the move, the heights of the crossings are 
\[h(a_l)\ge\ldots\ge h(a_1)\ge h(a)>h(b)\ge h(b_1)\ge\ldots\ge
h(b_m),\]
and let 
\begin{equation*}
b'_j\mapsto\varphi(b_j')=b_j\end{equation*} 
for all $j=1,\ldots,m$. Further\footnote{In this paper, each time we compute 
the holonomy of a Reidemeister II move, Proposition 3.5 
of \cite{en} applies. In other words, the rest of this section is mainly 
included only for completeness.}, write 
$\partial(a_1)=\sum B_1 b B_2 b B_3 b B_4b\ldots B_k b A$,
where $k\ge0$; 
$B_1,B_2,\ldots,B_k\in T(b_1,\ldots,b_m)$ are monomials; in the   
monomial $A\in T(b_1,\ldots,b_m,b,a)$, every $b$ factor is preceded by 
an $a$
factor (i.e., the $b$ factor right before $A$ is the
last $b$ before the first $a$; if there is no $a$ at all, then
it's just the last $b$ in the word); and the sum is taken over all
admissible discs with positive corner at $a_1$. Then let
\begin{equation*}
\begin{split}
a_1'\mapsto\varphi(a'_1)=a_1+\sum\Big(&
B_1aB_2bB_3bB_4b\ldots B_kbA\\
+&B_1vB_2aB_3bB_4b\ldots B_kbA\\
+&B_1vB_2vB_3aB_4b\ldots B_kbA\\
+&B_1vB_2vB_3vB_4a\ldots B_kbA\\
+&\ldots\\
+&B_1vB_2vB_3vB_4v\ldots B_kaA\Big).
\end{split}\end{equation*}
Finally, if the images $\varphi(a'_1),\ldots,\varphi(a'_{i-1})$ are 
already known and 
$\partial(a_i)=\sum B_1 b B_2 b B_3 b B_4b\ldots B_k b A$, where
$B_1,\ldots,B_k\in
T(b_1,\ldots,b_m,a_1,\ldots,a_{i-1})$, and $A$ doesn't contain any
copy of $b$ which is not preceded by a copy of $a$, then we let 
\begin{equation*}
\begin{split}
a_i'\mapsto\varphi(a_i')=a_i+\sum\Big(&
\bar{B_1}aB_2bB_3bB_4b\ldots B_kbA\\
+&\bar{B_1}v\bar{B_2}aB_3bB_4b\ldots B_kbA\\
+&\bar{B_1}v\bar{B_2}v\bar{B_3}aB_4b\ldots B_kbA\\
+&\bar{B_1}v\bar{B_2}v\bar{B_3}v\bar{B_4}a\ldots B_kbA\\
+&\ldots\\
+&\bar{B_1}v\bar{B_2}v\bar{B_3}v\bar{B_4}v\ldots\bar{B_k}aA\Big),
\end{split}\end{equation*}
where $\bar{B_1},\ldots,\bar{B_k}$ are obtained from $B_1,\ldots,B_k$
by replacing each symbol $a_j\in\{\,a_1,\ldots,a_{i-1}\,\}$ by the 
corresponding polynomial $\varphi(a'_j)$.

\section{Positive braid closures}\label{sec:BC}
Let us consider a positive braid $\beta$ on $q$ strands, as on Figure 
\ref{fig:fonat}. Label the left and right endpoints of the strands from 
top to bottom with the first $q$ whole numbers. The pair of left and right 
labels on each strand takes the form $(i,\sigma(i))$, where $\sigma$ is the 
\emph{underlying permutation} of $\beta$. Further, label the crossings 
of the braid with a pair of numbers, the first one the left label $i$ of the 
overcrossing strand and the second one the right label $j$ of the 
undercrossing 
one. (Note that not all pairs of numbers between $1$ and $q$ occur as labels 
of crossings; 
for example, $(i,\sigma(i))$ never does for any $i$.)

\begin{figure}
\psfrag{a1}{\Large $a_1$}
\psfrag{a2}{\Large $a_2$}
\psfrag{a3}{\Large $a_3$}
\psfrag{a4}{\Large $a_4$}
\psfrag{t1}{\Large $T_1$}
\psfrag{t2}{\Large $T_2$}
\psfrag{t3}{\Large $T_3$}
\psfrag{t4}{\Large $T_4$}
\psfrag{u1}{\Large $U_1$}
\psfrag{u2}{\Large $U_2$}
\psfrag{u3}{\Large $U_3$}
\psfrag{u4}{\Large $U_4$}
\psfrag{x1}{\Large $\alpha_1$}
\psfrag{x2}{\Large $\alpha_2$}
\psfrag{x3}{\Large $\alpha_3$}
\psfrag{x4}{\Large $\alpha_4$}
\psfrag{p1}{\Large $p_1$}
\psfrag{p2}{\Large $p_2$}
\psfrag{p3}{\Large $p_3$}
\psfrag{p4}{\Large $p_4$}
\resizebox{\linewidth}{521pt}{\includegraphics{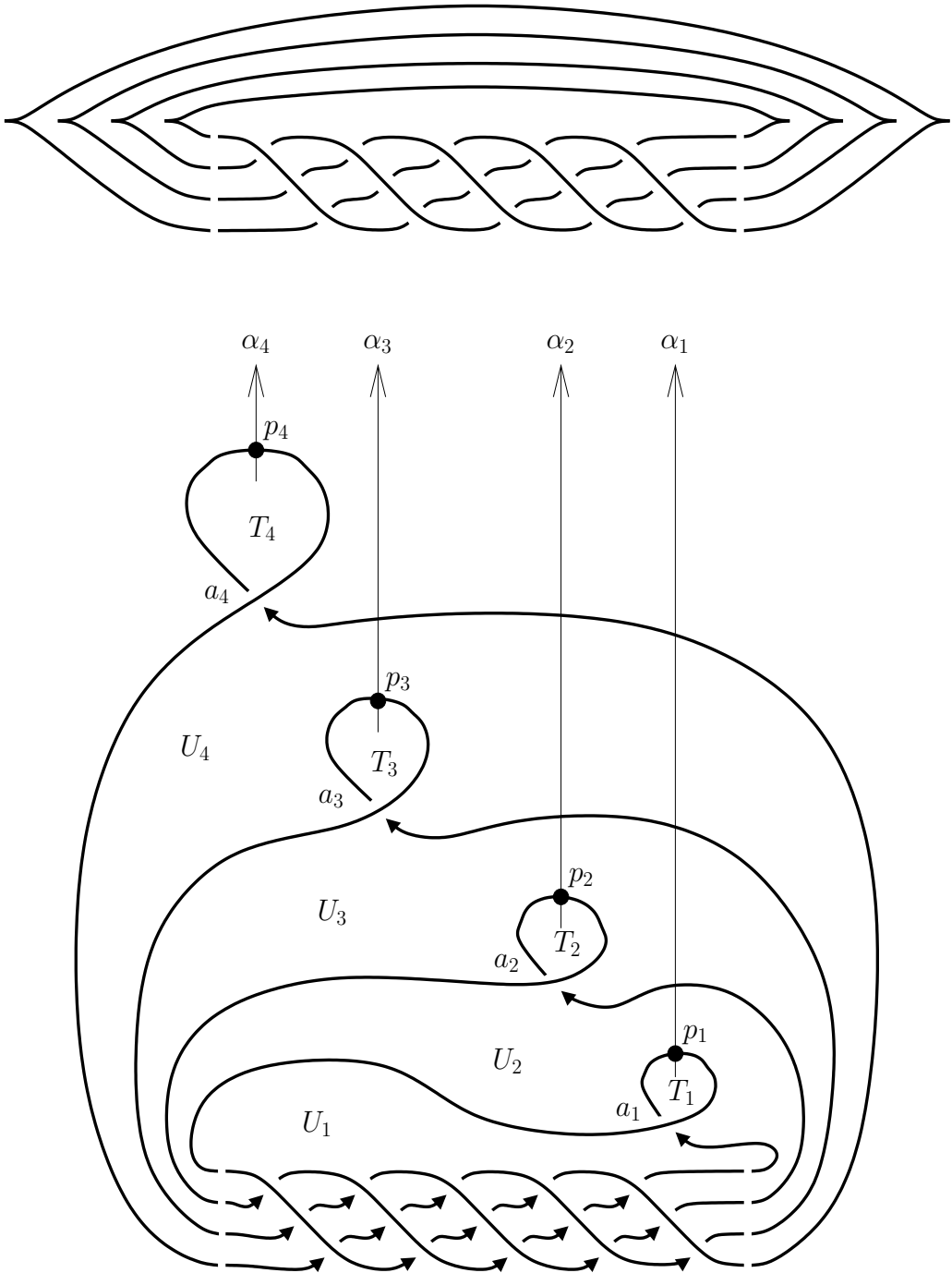}}
\caption{Front and Lagrangian diagrams of the closure of a positive 
braid}\label{fig:pic}
\end{figure}

\begin{Def}\label{def:lezaras}
A positive braid $\beta$ defines the front diagram of an oriented 
Legendrian link as on the upper half of Figure \ref{fig:pic} 
(orient each strand of $\beta$ from left to right). We call 
this construction the \emph{Legendrian closure} of $\beta$ and denote it by 
$L_\beta$. 
\end{Def}


Note that the types of 
crossings are correctly determined by the slopes of the branches that meet 
there exactly because the braid is positive. 
Applying Ng's \cite{computable} construction of 
resolution to the front diagram, we obtain the Lagrangian diagram 
$\gamma_\beta$ on the lower half of Figure \ref{fig:pic}. 


In section 2 
of \cite{en}, we omitted the modifications needed 
to define the invariants (such as rotation and Thurston--Bennequin numbers and 
contact homology itself) for oriented multi-component links. We include an 
informal rundown here (see \cite[section 2.5]{computable} for more). The 
rotation number is the sum of the rotations of the 
components. The Thurston--Bennequin number is still the writhe of the 
Lagrangian projection; note that now, as opposed to the case of a single 
component, the orientation matters in its definition. There is not a 
single distinguished grading of a link DGA, rather a family of so 
called admissible gradings. 
We consider those introduced in \cite{computable} and not the larger 
class of gradings described in \cite[section 9.1]{chek}.
In the $\Z_2$--coefficient theory, each of these is defined modulo the 
greatest common divisor of the 
Maslov numbers of the components (recall that the Maslov number is twice 
the rotation number of a knot). 
Self-crossings of individual components 
have the same well-defined index in any admissible grading. We will refer to 
such generators as \emph{proper crossings}. Other crossings' indices have the 
same parity in each admissible grading. The parity of the index of any 
crossing coincides with the sign in classical knot theory.

The differential $\partial$ is of index 
$-1$ with respect to each admissible grading. Contact homology is a 
well-defined invariant in the sense that if two diagrams are Legendrian 
isotopic, then there is a one-to-one correspondence between their sets 
of admissible gradings and the corresponding contact homologies are 
isomorphic as graded algebras.

It is easy to calculate 
that the Legendrian closure of any positive braid has rotation number $r=0$ 
(in fact, every component has Maslov number $0$). The Thurston--Bennequin 
number is $tb(L_\beta)=(\text{word length of }\beta)-q$. 
The positive Legendrian $(p,q)$ 
torus knot obtained as a special case 
is the one with maximal Thurston--Bennequin number $p(q-1)-q$ (see 
\cite{EH1} for a classification of Legendrian torus knots). The $(3,2)$ 
torus knot considered in \cite{en} arises from this construction, too. In 
fact, the unique Legendrian unknot with maximal Thurston--Bennequin 
number $tb=-1$ \cite{benn,unknot} is the Legendrian closure of the trivial 
braid on a single strand.

\begin{figure}
\psfrag{bij1}{\LARGE $b_{i,j,1}$}
\psfrag{bij2}{\LARGE $b_{i,j,2}$}
\psfrag{bij3}{\LARGE $b_{i,j,3}$}
\psfrag{bpq1}{\LARGE $b_{p,r,1}$}
\psfrag{bpq2}{\LARGE $b_{p,r,2}$}
\psfrag{i}{\LARGE $i$}
\psfrag{j}{\LARGE $j$}
\psfrag{p}{\LARGE $p$}
\psfrag{q}{\LARGE $r$}
\resizebox{\linewidth}{!}{\includegraphics{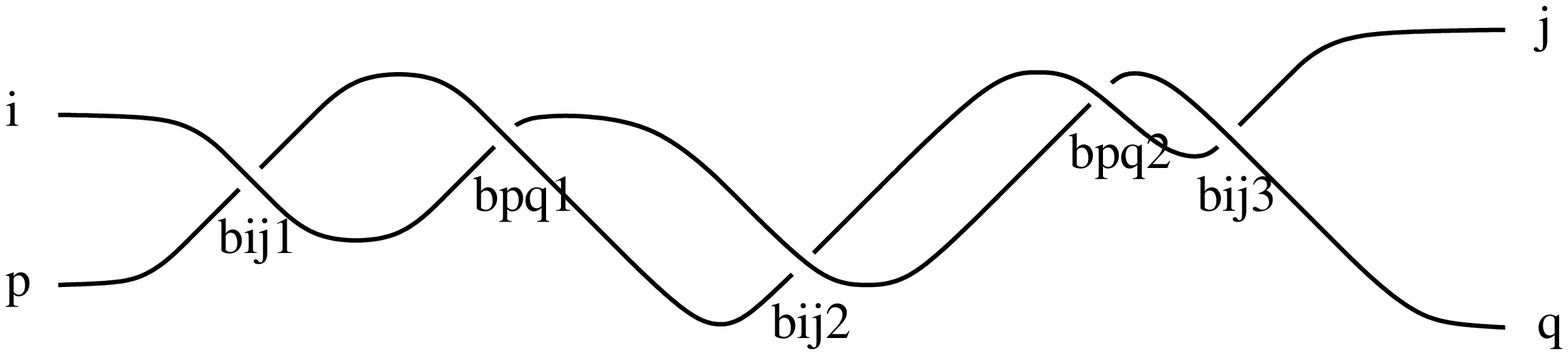}}
\caption{Labeling the crossings of a braid}\label{fig:indexek}
\end{figure}

In the Chekanov--Eliashberg DGA of $\gamma_\beta$, the grading 
(any admissible grading) is integer-valued.
The crossings $a_m$ 
($m=1,\ldots,q$) have index $1$. 
The rest of the generators are the crossings of $\beta$ 
and they have index $0$. In the multi-component case we should say instead 
that 
the grading that assigns the index $0$ to each is 
admissible; from now on, we will always work with this grading. 
These generators will be labeled by $b_{i,j,t}$, where 
the integers $1\le i,j\le q$ are the ones defined at the beginning 
of the section. The third label $t$ is used to distinguish between 
multiple intersections of strands, enumerating them from left to right, as 
shown on Figure \ref{fig:indexek}.

Since $\partial$ lowers the index by $1$, 
\[\partial(b_{i,j,t})=0\text{ for all }i,j=1,\ldots,q\text{ and }t.\]
In fact, there are no admissible discs whose positive corner is an index 
$0$ crossing. The boundaries $\partial(a_m)$, $m=1,\ldots,q$, are 
polynomials in the non-commuting variables $b_{i,j,t}$ only. Our next 
goal is to compute these polynomials.

The first observation is that for all $m=1,\ldots,q$, there is an 
admissible disc covering the teardrop-shaped region $T_m$ right above $a_m$ 
once. We'll call it the \emph{$m$'th trivial disc}. It contributes $1$ to 
$\partial(a_m)$. The $m$'th trivial disc is the only one that turns at the 
upward-facing positive quadrant at $a_m$. The rest of the contributions to 
$\partial(a_m)$ come from discs that turn at the positive quadrant facing 
down.

Suppose $f\colon\Pi_k\to\R^2_{xy}$ is an admissible immersion with respect 
to the projection $\gamma=\gamma_\beta$, with 
positive corner $f(x_0^k)=a_m$ and 
so that it is different from the $m$'th trivial disc. Fix 
points $p_1,\ldots,p_q$ on $\gamma$ as shown on Figure \ref{fig:pic}. 

\begin{lemma}\label{lem:kicsi}
The curve $f(\partial\Pi_k)$ doesn't pass through the points 
$p_m,\ldots,p_q$.
\end{lemma}

\begin{proof}
In fact, if we denote the region of the 
complement directly under $a_i$ by $U_i$, then 
$f(\Pi_k)$ is disjoint from $U_{m+1},\ldots,U_q$, as well as from 
$T_m,\ldots,T_q$. To see this, first shrink $\gamma$ so that all areas and 
heights are smaller than $\varepsilon^{q-m}$, where $\varepsilon>0$ is to be 
chosen later. Then add ``bulges'' of equal area $1$ to $T_q$ and $U_q$ on 
the outside of the diagram (see Remark 2.8
of \cite{en}). Next, add bulges of area $\varepsilon$ to $T_{q-1}$ and 
$U_{q-1}$ on their sides facing $U_q$. Continue all the way until adding 
bulges of area $\varepsilon^{q-m-1}$ to $T_{m+1}$ and $U_{m+1}$ on their 
sides facing $U_{m+2}$. The result is a Lagrangian projection of the 
same link in which, 
after the choice of a small enough $\varepsilon$, the 
claim is obvious by Lemma 2.11
of \cite{en}, except for the case of $T_m$. But that follows from the easy 
observation that if any admissible disc passes through $p_m$ so that locally, 
$f(\Pi_k)$ faces downward, then $f$ is the $m$'th trivial disc.
\end{proof}

Recall that $\gamma$ is oriented as shown on Figure \ref{fig:pic}.
We have already noted that apart from its positive corner which is 
of index $1$, $f$ turns only at certain index $0$ crossings 
$b_{i,j,t}$. By Lemma 2.10
of \cite{en}, this implies that $f$ is compatible with the 
orientation of $\gamma$. Therefore, as the orientation of 
$f(\partial\Pi_k)$ 
agrees with that of $\gamma$ near the positive corner $a_m$, these 
orientations agree at all points of $\partial\Pi_k$. 
This last observation for example implies that $f$ can't turn at the 
downward facing negative quadrant at any of the $b_{i,j,t}$'s. 

We may summarize our findings about $f$ so far as follows. The curve 
$f(\partial\Pi_k)$ 
starts at $a_m$, follows $\gamma$ until it reaches the braid at the left 
endpoint labeled $m$ (while $f$ extends this map toward $U_m$). 
Then it travels through the braid, 
possibly turning (left) at several crossings $b_{i,j,t}$ but always 
heading 
to the right, until it reaches a right endpoint labeled $i_1\le m$. Then 
$f(\partial\Pi_k)$ climbs to $a_{i_1}$ and, unless $i_1=m$, beyond 
$a_{i_1}$ (note that $f$ can't have a positive corner at $a_{i_1}$) 
to $p_{i_1}$, at which point the extension $f$ is toward $U_{i_1+1}$. Then 
$f(\partial\Pi_k)$ descends back to $a_{i_1}$ so that it doesn't turn 
at the negative corner. So the process repeats with $(m,i_1)$ 
replaced with $(i_1,i_2)$ (where $i_2\le m$), 
and so on until $i_{c+1}=m$ for some $c$. 

\begin{Def}\label{def:sorozat}
A finite sequence of positive integers is called \emph{admissible} if for 
all $s\ge0$, between any two appearances of $s$ in the sequence there is a 
number greater than $s$ which appears between them. For $n\ge1$, let us 
denote by $D_n$ 
the set of all admissible sequences that are composed of the numbers 
$1,2,\ldots,n-1$. 
\end{Def}

For example, $D_1=\{\:\varnothing\:\}$, 
$D_2=\{\:\varnothing,\{\,1\,\}\:\}$, and 
\[D_3=\{\:\varnothing,
\{\,1\,\},\{\,2\,\},\{\,1,2\,\},\{\,2,1\,\},\{\,1,2,1\,\}\:\}.\] 
By induction on $n$ 
and observing the position of the unique maximal term in the sequence, it 
is easy to prove that $|D_n|=|D_{n-1}|^2+|D_{n-1}|$.

\begin{all}\label{pro:blank}
For the admissible disc $f$ as above, the 
sequence of intermediate labels 
$\{\,i_1,i_2,\ldots,i_c\,\}$ is an element of $D_m$.
\end{all}

\begin{proof}
There is an intuitive reason for the claim: if the index $i$ was repeated 
without the boundary of the disc climbing higher between the two 
occurrences, then the disk would pinch off at $a_i$. We shall give a more 
rigorous proof using Blank's theorem 
(see \cite{blank}, or the review 
in \cite{franc}, from where we'll borrow our terminology). 

Let $\tilde f$ be a small generic perturbation of $f$. 
This has the corners rounded so that $\tilde f\big|_{\partial\Pi_k}$ is an 
immersion; moreover, let us note that 
because $f(\partial\Pi_k)$ may cover parts of $\gamma$ more than once, 
the complement of $\tilde f(\partial\Pi_k)$ typically contains more 
regions than that of $f(\partial\Pi_k)$. 
We are going to apply Blank's theorem to 
$\tilde f\big|_{\partial\Pi_k}$. 
We extend rays from each bounded component of 
$\R^2\setminus \tilde f(\partial\Pi_k)$ to infinity. 
If the region is not one of those obtained from the $T_i$ 
(let us denote these by $T'_i$), or one of the 
small ones resulting from perturbation
close to the arc bounding $T_i$, this can be done so that $\gamma$ 
only intersects the ray in a positive manner (i.e., if the ray is 
oriented toward $\infty$, then $\gamma$ crosses from the right side to 
the left side). By our observation on the compatibility of $f$, this 
implies that $\tilde f(\partial\Pi_k)$ also intersects those rays positively. 
Let us draw rays starting from the regions $T'_1,\ldots,T'_q$ as on Figure 
\ref{fig:pic}. Label all rays, in particular label these last $q$ with 
the symbols $\alpha_1,\ldots,\alpha_q$.

The \emph{Blank word} of the disc $\tilde f$ is obtained by tracing 
$\tilde f(\partial\Pi_k)$, starting from, say, $a_m$, and writing down the 
labels of 
the rays we meet with exponents $\pm1$ 
according to positive and negative 
intersections (we'll call these 
\emph{positive} and \emph{negative symbols} or \emph{letters}). A 
\emph{grouping} of the Blank word is a set of 
properly nested disjoint unordered pairs of 
the form $\{\,\alpha,\alpha^{-1}\,\}$ so that each negative symbol 
is part of exactly one pair.
Blank's theorem states that the set of groupings of the Blank word is in 
a one-to-one correspondence with non-equivalent extensions of the 
immersion $\tilde f\big|_{\partial\Pi_k}$ to immersions of $\Pi_k$. 

As $\tilde f$ is such an extension, a grouping exists and we claim that this 
implies the Proposition. We delete labels different from the $\alpha_s$ 
from the Blank word of $\tilde f$ and concentrate on 
the remaining word, which also inherits a grouping (note that most of the 
deleted labels only appeared positively anyway, except for certain ones 
that belong to some regions that were the result of perturbation). 
We'll call this the \emph{Blank word of $f$} and 
denote it by $W_f$. Note that if we only keep the negative letters 
$\alpha_s^{-1}$ from $W_f$, then the sequence of their indices is exactly 
$\{\,i_1,i_2,\ldots,i_c\,\}$. In fact, $W_f$ is decomposed into segments 
$S_1,\ldots,S_c$ ended by these negative symbols, and a 
final segment $S$:
\begin{multline*}
W_f=
\overbrace{\alpha_1\alpha_2\ldots\alpha_{i_1-1}\alpha_{i_1}^{-1}}^{S_1}
\cdot
\overbrace{\alpha_1\alpha_2\ldots\alpha_{i_2-1}\alpha_{i_2}^{-1}}^{S_2}
\cdot\ldots\\
\cdot\underbrace{\alpha_1\alpha_2\ldots\alpha_{i_c-1}\alpha_{i_c}^{-1}}
_{S_c}\cdot
\underbrace{\alpha_1\alpha_2\ldots\alpha_{m-1}}_S.
\end{multline*}

We'll prove the following statements by induction on $j$: 
\begin{enumerate}[(1)]
\item\label{rend} 
Each copy of $\alpha^{-1}_{m-j}$ is paired in the grouping with a copy 
of $\alpha_{m-j}$ which is located to the right of it.
\item\label{benn} 
None of the pairs $\{\,\alpha^{-1}_{m-j},\alpha_{m-j}\,\}$ is nested in 
a pair $\{\,\alpha^{-1}_{m-l},\alpha_{m-l}\,\}$ for any $l>j$.
\item\label{jo} 
Any two copies of $\alpha^{-1}_{m-j}$ are separated in $W_f$ by a copy 
of $\alpha^{-1}_{m-n}$ for some $n<j$.
\end{enumerate}
The last statement of course directly implies the Proposition.

The last letter of $W_f$ is $\alpha_{m-1}$ and this is the only positive 
occurrence of this letter. 
Hence \eqref{rend} and \eqref{benn} are obvious for $j=1$. 
Also, since a second one wouldn't find a pair, 
there can be at most one copy of $\alpha^{-1}_{m-1}$ in $W_f$, thus 
\eqref{jo} is vacuously true for $j=1$. 

Assume that the statements hold for all $j'=1,\ldots,j-1$. To prove 
\eqref{rend}, assume that a certain copy $\alpha^{-1}$ of 
$\alpha^{-1}_{m-j}$ forms the 
pair $P$ with a copy $\alpha$ of $\alpha_{m-j}$, which is to the 
left of it. Then, as 
$\alpha$ can't be part of the final segment $S$, it belongs to some 
$S_b$ with last letter $\alpha^{-1}_{i_b}$, where $i_b>m-j$, i.e.\ 
$m-i_b<j$. This copy of $\alpha^{-1}_{i_b}=\alpha^{-1}_{m-(m-i_b)}$ is 
to the left of $\alpha^{-1}$, hence it is part of a pair which is nested 
inside $P$, which contradicts the hypothesis \eqref{benn} for $j'=m-i_b$.

To prove \eqref{benn}, suppose that a certain pair 
$\{\,\alpha^{-1}_{m-j},\alpha_{m-j}\,\}$ is nested in $P=\{\,\alpha, 
\alpha^{-1}\,\}$ where $\alpha$ has an index less than $m-j$. This copy 
of $\alpha_{m-j}$ can't be part of $S$ because it has a symbol with lower 
index (namely, $\alpha$ or $\alpha^{-1}$) to the right of it. It is clear 
then that the first negative letter after $\alpha_{m-j}$ has index 
$m-j'$, which is higher
than $m-j$ (i.e.\ $j'<j$), and it is still nested in $P$. Then so is the 
pair containing it, which contradicts \eqref{benn} for $j'$.

Finally for \eqref{jo}, assume the contrary again, namely that there are 
two copies of $\alpha^{-1}_{m-j}$ in $W_f$ that are not separated by any 
higher index negative symbol. Then the pair of the first 
$\alpha^{-1}_{m-j}$ can not lie between them either, because it would be 
part of some $S_b$ and then $\alpha^{-1}_{i_b}$ would separate. Thus, we 
have two pairs of the form $\{\,\alpha^{-1}_{m-j},\alpha_{m-j}\,\}$ nested 
in one another. But then the positive symbol of the inner pair would be 
followed by some $\alpha^{-1}_{i_b}$, which is part of a pair that is 
nested in the outer pair, and that contradicts \eqref{benn} for 
$j'=m-i_b$.
\end{proof}

\begin{Def}\label{def:B}
Let $1\le i,j\le q$. The element $B_{i,j}$ of the DGA of $\gamma_\beta$ 
is the sum of the following 
products. For each path composed of parts of the strands of the braid $\beta$ 
that connects the left endpoint labeled $i$ to the 
right endpoint labeled $j$ so that it only turns at quadrants facing up, 
take the product of the crossings from left to right that it turns at.
\end{Def}

For example, $B_{i,j}$ contains the constant term $1$ if and only if 
$j=\sigma(i)$. We will need the following polynomials of the $B_{i,j}$:

\begin{Def}\label{def:C}
Let $q\ge i>j\ge1$ and let 
\[C_{i,j}=\sum_{\{\,i_1,\ldots,i_c,j\,\}\in D_i}
B_{i,i_1}B_{i_1,i_2}B_{i_2,i_3}\ldots B_{i_{c-1},i_c}B_{i_c,j}.\]
Similarly, for $1\le m\le q$, let 
\begin{equation}\label{eq:Cii}
C_{m,m}=\sum_{\{\,i_1,\ldots,i_c\,\}\in D_m}
B_{m,i_1}B_{i_1,i_2}B_{i_2,i_3}\ldots B_{i_{c-1},i_c}B_{i_c,m}.
\end{equation}
Finally, for any $i$ and $j$, let
\begin{equation*}
M_{i,j}=\sum_{\{\,i_1,\ldots,i_c\,\}\in D_{\min\{i,j\}}}
B_{i,i_1}B_{i_1,i_2}B_{i_2,i_3}\ldots B_{i_{c-1},i_c}B_{i_c,j}.
\end{equation*}
\end{Def}

In particular, for each summand in $C_{i,j}$, $j$ is the last element of 
the admissible sequence but it may occur elsewhere, too. For example, 
$C_{1,1}=B_{1,1}$, $C_{2,1}=B_{2,1}$, and 
$C_{3,1}=B_{3,1}+B_{3,2}B_{2,1}+B_{3,1}B_{1,2}B_{2,1}$. 
Note also that $M_{1,j}=B_{1,j}$, $M_{i,1}=B_{i,1}$, $M_{m,m}=C_{m,m}$ and 
$M_{i,i-1}=C_{i,i-1}$, whenever these expressions are defined.


\begin{tetel}\label{thm:relacio}
$\partial(a_m)=1+C_{m,m}$. Consequently, the index $0$ part 
$H_0(L_\beta)$ of the contact 
homology $H(L_\beta)$ has a presentation where the generators are the 
crossings of $\beta$ and the relations are $C_{m,m}=1$ for 
$m=1,\ldots,q$. 
\end{tetel}

\begin{proof}
The $1$ in the formula comes from the $m$'th trivial disc. We claim that 
the rest of the 
contributions add up to $C_{m,m}$. From Proposition \ref{pro:blank} and the 
paragraph preceding Definition \ref{def:sorozat}, we know that there can't 
be any such terms other than the ones included in $C_{m,m}$. To see that all 
such monomials actually arise from admissible discs, we just need to find 
those discs. This can either be done by an inductive construction on $c$ (for 
the inductive step, remove the smallest number from the sequence), or by 
applying Blank's theorem.
\end{proof}

It is not clear whether $H(L_\beta)$ contains any non-zero higher index 
part at all (except for the case of the unknot, when the single index $1$ 
crossing is a non-nullhomologous cycle). This is mainly why we only 
work with $H_0(L_\beta)$ in this paper. This difficulty in handling 
contact homology also underlines the importance of the augmentation that 
we construct in the next section.

\section{Augmentations of braid closures}\label{sec:aug}
\begin{Def}
An \emph{augmentation} of the Lagrangian diagram $\gamma$ of a Legendrian 
link $L$ is a subset $X$ of its crossings with the following properties.
\begin{itemize}
\item All elements of $X$ are proper crossings of $L$ 
(i.e., intersections of different components are not allowed in $X$). 
\item The index of each element of $X$ in any admissible grading is $0$ 
(in fact, this requirement implies the previous one).
\item For each generator $a$, the number of admissible discs with positive 
corner $a$ and all negative corners in $X$ is even.
\end{itemize}
\end{Def}

The last requirement implies that the evaluation homomorphism (which is 
defined on the link DGA, and which is also called an augmentation) 
$\varepsilon_X\colon\mathscr A\to\Z_2$ that sends elements of $X$ 
to $1$ and other generators to $0$, gives rise to an algebra homomorphism 
$(\varepsilon_X)_*\colon H(L)\to\Z_2$.


\begin{pelda}\label{ex:aug}
Consider the right-handed Legendrian trefoil knot diagram on Figure 4
of \cite{en}. We claim that the set $\{\,b_3\,\}$ is an augmentation. 
Indeed, the only two nonzero differentials (see Example 2.14
of \cite{en} or the previous section) are $\partial(a_1)=1+b_1+b_3+b_1b_2b_3$ 
and $\partial(a_2)=b_2+b_2b_3+b_1b_2+b_2b_3b_1b_2$, and even these vanish 
after mapping $b_1$ and $b_2$ to $0$ and $b_3$ to $1$.
\end{pelda}


\begin{megj}\label{rem:fuchs}
Let us return to the front diagram in the upper half of Figure \ref{fig:pic}.
Such a diagram always has an \emph{admissible decomposition} (or \emph{ruling})
in the sense of \cite{chek2}: the only two values
of the Maslov potential (even though it's $\Z$--valued) 
are $1$ on the upper strands and $0$ on the strands of the original braid 
$\beta$, thus all crossings are Maslov, and 
we may declare all of them switching (in the multi-component case, 
consider only proper crossings). This gives 
rise to a decomposition where the discs are nested in one another, so it's 
admissible.

The existence of a ruling implies (see \cite{chek2}) 
that $L_\beta$ is not a stabilized link type 
for any positive braid $\beta$.
Also, by a theorem of Fuchs \cite{fuchs}, it implies that the diagram has an 
augmentation (note this is by no means unique). In his proof, Fuchs 
constructs an augmentation of a diagram which is equivalent 
to the original but has a lot more crossings. This means that the original 
diagram can also be augmented: the pull-back of an augmentation by a DGA 
morphism, like the ones listed in section \ref{sec:z2}, is again an 
augmentation. 
\end{megj}

\begin{figure}
\psfrag{1}{\LARGE $1$}
\psfrag{2}{\LARGE $2$}
\psfrag{3}{\LARGE $3$}
\psfrag{4}{\LARGE $4$}
\psfrag{5}{\LARGE $5$}
\psfrag{6}{\LARGE $6$} 
\psfrag{7}{\LARGE $7$} 
\psfrag{8}{\LARGE $8$} 
\psfrag{a}{\large $1$}
\psfrag{b}{\large $2$}
\psfrag{c}{\large $3$}
\psfrag{d}{\large $4$}
\psfrag{e}{\large $5$}
\psfrag{f}{\large $7$}
\resizebox{\linewidth}{!}{\includegraphics{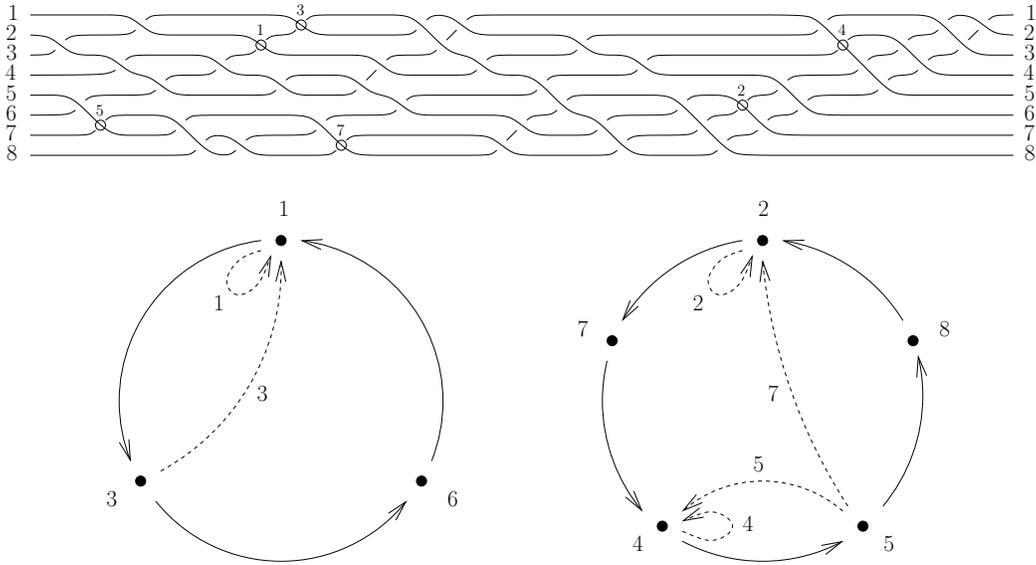}}
\caption{Constructing an augmentation for the closure of a 
braid}\label{fig:fonat}
\end{figure}

In the case of a braid closure, it would be 
impractical to pull back Fuchs' augmentation to the original diagram. Instead, 
we'll start from scratch and construct an augmentation of the Legendrian 
closure $\gamma_\beta$ of an arbitrary positive braid $\beta$. 
When selecting crossings into $X$, it will suffice to work with the 
braid itself, as illustrated on Figure \ref{fig:fonat}, for the index $0$ 
proper crossings of $\gamma_\beta$ are all crossings of $\beta$. 
(We will call these the \emph{proper crossings of the braid $\beta$}.)

First, we associate an oriented graph to an arbitrary 
permutation $\sigma\in S_q$. Let $o$ be one of the cycles of 
$\sigma$. Let us write the elements of $o$ around the perimeter of a 
circle in the 
cyclic order suggested by $\sigma$, directing an edge from $s$ to 
$\sigma(s)$ for all $s$. (See Figure \ref{fig:fonat} for an example. The 
permutation on the diagram is the one underlying the braid). 
If $p$ is a non-maximal element of 
$o$, then follow the cycle in the forward direction starting from $p$ 
until it hits the first number in $o$ which is bigger then $p$. Let 
$p_+$ be the number immediately before that. In particular, $p_+\le 
p<\sigma(p_+)$. Do the same in the backward direction, resulting in the 
number $p_-$ such that $p_-\le p<\sigma^{-1}(p_-)$. 
(For example, on Figure \ref{fig:fonat}, $5$ is an element of a $5$--cycle 
$o$, with $5_+=5$ and $5_-=4$.) Then, connect $p_+$ to 
$p_-$ by a directed chord of the circle and label the chord by $p$. 
If $p_+=p_-=p$, then instead, we attach a loop edge (labeled $p$) 
at $p$ to the graph. This will also be called a chord. Draw this loop 
edge inside the circle, right next to the perimeter, 
on the side of $p$ where the smaller of its two neighbors lies. 

\begin{Def}
If $p$ is non-maximal in its cycle $o$, then the oriented loop $\Gamma_p$ 
that starts from $p$, goes along the circle of $o$ to $p_+$, then 
goes to $p_-$ on 
the chord labeled $p$, then follows the circle again back to $p$ will be 
called the \emph{loop of $p$}. If $\bar p$ is the largest number in $o$, 
then let the \emph{loop of $\bar{p}$} be the loop $\Gamma_{\bar p}$ that 
travels around the original circle once.
\end{Def}

\begin{lemma}\label{lem:diszj}
If $p<r$, then 
\begin{enumerate}[(a)]
\item\label{ize} $\Gamma_p$ doesn't contain $r$
\item\label{bize} the discs bounded by $\Gamma_p$ and $\Gamma_r$ are either 
disjoint or the latter contains the former.
\end{enumerate}
In particular, the $|o|-1$ (oriented) chords obtained in the construction are 
pairwise disjoint. They can't be parallel to the original edges and 
they differ from each other as well.
\end{lemma}

\begin{proof}
Chords differ from edges because it's impossible that each end of an edge 
be smaller than the other end (a chord and an edge can be parallel in the 
non-oriented sense, as on Figure \ref{fig:fonat}). The other statements 
follow from \eqref{ize} and \eqref{bize}. Statement \eqref{ize} is obvious 
from the construction, and so is \eqref{bize} if $r$ is maximal in $o$. 

The directed 
arc of the circle stretching from $p_-$ to $p_+$ only contains numbers less 
than $r$. If this arc is disjoint from $\Gamma_r$, then of course so is 
$\Gamma_p$. Otherwise, the whole arc must be contained in $\Gamma_r$. From 
this, statement \eqref{bize} is clear, except if the chord labeled $p$ is a 
loop edge and either $p=r_+$ or $p=r_-$. In the first case, $\sigma(p)>r$, but 
$\sigma^{-1}(p)\le r$, so \eqref{bize} holds by the construction of the 
loop edge. The other case is handled analogously.
\end{proof}

\begin{Def}
Let $\Gamma_\sigma$ be the disjoint union of 
the graphs constructed above over all cycles $o$ of $\sigma\in S_q$. 
This oriented planar graph, with vertices labeled
by the numbers $1,\ldots,q$, is called the \emph{augmented graph of the 
permutation $\sigma$}.
\end{Def}

\begin{lemma}\label{lem:graf-szint}
For all $p\in\{\,1,\ldots,q\,\}$, the loop of $p$ is the 
unique directed loop in $\Gamma_\sigma$ 
starting and ending at $p$ so that the sequence of the vertices (apart from 
$p$) visited by it is in $D_p$.
\end{lemma}

\begin{proof}
It is enough to prove the statement 
for a connected component associated to a cycle $o$. 
The loop $\Gamma_p$ 
has the required property for all $p$, because all the vertices visited 
by it are less than $p$ and no repetition occurs. 

We have to rule out the existence of other loops. 
Since $\sigma(p_+)$ and $\sigma^{-1}(p_-)$ are larger 
than $p$, by the disjointness statement in Lemma \ref{lem:diszj}, 
all loops in question are trapped 
in the disc bounded by $\Gamma_p$ (and this is obviously true 
when $p$ is maximal in $o$). Then by statement \eqref{bize} of 
Lemma \ref{lem:diszj}, apart from edges of $\Gamma_p$, they may 
only contain chords labeled by numbers less than $p$. 
Now, if the chord labeled by $s$ occurred in the loop and $s<p$ was the 
smallest such number, then because $p$ is not on $\Gamma_s$, 
the sequence of vertices on $\Gamma_s$ would appear as 
a subsequence of the original vertex sequence of the loop. In that case, 
$s$ would be 
repeated in the sequence without the two occurrences separated by a 
larger number, which is a contradiction.
\end{proof}

The following ``edge reversal lemma,'' which we will need in section 
\ref{sec:moreaug}, can be proven very similarly.

\begin{lemma}\label{lem:masodik}
For any $1\le p<r\le q$ so that $\sigma(r)=p$, if $p$ is the second 
largest vertex (after $r$) along $\Gamma_r$, then there is a unique 
oriented path in $\Gamma_\sigma$ from $p$ to $r$ so that the sequence of 
the intermediate vertices is in $D_p$. Otherwise, there is no such path. 
\end{lemma}

\begin{Def}
Let $Y$ be a set of crossings of the positive braid $\beta$. By the 
\emph{graph realized by 
$Y$} we mean the oriented graph with vertices $1,\ldots,q$ so that a 
directed edge 
connects $i$ to $j$ if and only if $B_{i,j}\big|_Y=1$. Here, $B_{i,j}$ is 
as in 
Definition \ref{def:B} and by $B_{i,j}\big|_Y$ we mean the element of 
$\Z_2$ obtained by substituting $1$ for elements of $Y$ and $0$ for other 
generators in $B_{i,j}$.
\end{Def}

\begin{lemma}\label{lem:grafok}
Let $\sigma$ be the underlying permutation of $\beta$ and $Y$ a set of 
proper crossings of $\beta$. If the graph realized by $Y$ agrees with the 
augmented graph $\Gamma_\sigma$ 
of $\sigma$, then $Y$ is an augmentation of the Legendrian 
closure of $\beta$.
\end{lemma}

\begin{proof}
Assume the two graphs do agree. Then by Lemma \ref{lem:graf-szint}, exactly 
one of the summands of 
$C_{m,m}$ (see equation \eqref{eq:Cii}) contributes $1$ to 
the sum $C_{m,m}\big|_Y$, 
namely the one that belongs to the sequence of vertices on 
$\Gamma_m$. Therefore by Theorem \ref{thm:relacio}, $Y$ is an 
augmentation.
\end{proof}

Next, based on $\Gamma_\sigma$, we construct a candidate $X$, and then we 
will use Lemma \ref{lem:grafok} to prove that it's an augmentation.
Loosely speaking, the edges of $\Gamma_\sigma$ 
connecting $s$ to $\sigma(s)$ are always realized, even by the 
empty set. To realize the chord labeled $p$, we'll select the crossing 
$b_{p_+,p_-,1}$ into $X$. We can do this because it always exists: 
$p_+\le 
p<\sigma^{-1}(p_-)$ (for left labels) and $p_-\le p<\sigma(p_+)$ (for 
right labels), therefore the strand connecting $p_+$ to $\sigma(p_+)$ 
always meets the strand connecting $\sigma^{-1}(p_-)$ to $p_-$. If there 
were more than one points with the first two labels $p_+,p_-$, we could have 
selected any of them\footnote{For instance, the crossing $b_3$ of Example 
\ref{ex:aug} is denoted by $b_{1,1,2}$ in the general labeling system.}; 
we used the third label $1$ for concreteness 
and for ease in the proof of Theorem \ref{thm:osztja}. On 
Figure \ref{fig:fonat}, we marked the selected crossings and labeled them 
with the label of the chord that they realize.

\begin{all}\label{pro:aug}
The set 
\[X=\{\:b_{p_+,p_-,1}\mid p\in\{\,1,\ldots,q\,\}\text{ is not a 
maximal element of a cycle of }\sigma\:\}\]
is an augmentation of the Legendrian closure of the positive braid $\beta$ 
with underlying permutation $\sigma$.
\end{all}

In particular, for a pure braid $\beta$, the empty set is an augmentation. 
In other words, the DGA of the Legendrian closure of a pure braid is 
augmented, i.e.\ the boundary of each generator is a polynomial without 
a constant term. If $\beta$ is not pure, then $X\ne\varnothing$, hence 
$\varepsilon_X\ne 0$, and it follows that $H(L_\beta)\ne 0$.

\begin{proof}
It is clear from the construction that all the crossings in $X$ are proper. 
By Lemma \ref{lem:grafok}, it suffices to prove that the graph $G$ 
realized by $X$ is the graph $\Gamma_\sigma$. For 
this, the chief claim is that no two points of $X$ are connected with a 
part of a strand so that it arrives at both points from above. In other 
words, the situation of Figure \ref{fig:baj} can not arise: there is no 
pair of numbers $p,r$ so that $\sigma(r_+)=p_-$. Indeed, then we'd have 
$r<\sigma(r_+)=p_-\le p$ and similarly, $p<\sigma^{-1}(p_-)=r_+\le r$, 
which would be a contradiction.

\begin{figure}
\psfrag{p+}{\Large $p_+$}
\psfrag{p-}{\Large $p_-$}
\psfrag{r+}{\Large $r_+$}
\psfrag{r-}{\Large $r_-$}
\psfrag{p}{\Large $p_+,p_-$}
\psfrag{r}{\Large $r_+,r_-$} 
\resizebox{\linewidth}{!}{\includegraphics{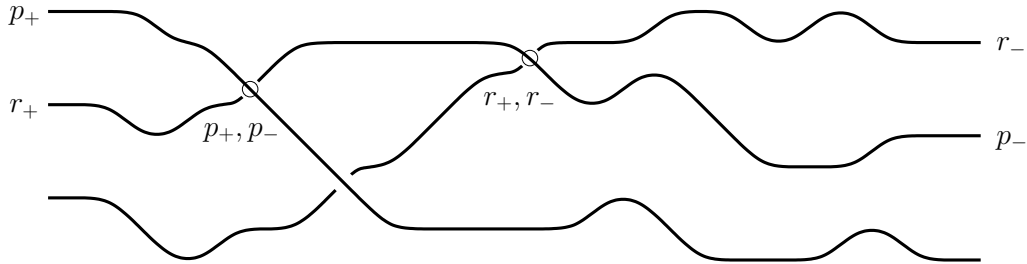}}
\caption{A situation that we have to rule out}\label{fig:baj}
\end{figure}

We know then that any path that is to contribute a non-zero summand to a 
certain $B_{i,j}\big|_X$ can have at most one corner, which of course 
has to be in $X$. So in the discussion before the Proposition we exhausted 
all such contributions: paths with no corners are responsible for the 
edges and paths with one corner are responsible for the chords of 
$\Gamma_\sigma$. So indeed, $G=\Gamma_\sigma$. 
\end{proof}

\begin{megj}
Let $\beta$ be the standard positive braid whose Legendrian closure is a 
positive $(p,q)$ torus link (the braid used to produce Figure 
\ref{fig:pic} is an example with $p=5$ and $q=4$).
If we apply our construction to it, we find 
an interesting connection of the resulting augmentation to the Euclidean 
algorithm. We mention this here 
without proof; we will only need a small part of the statement in section 
\ref{sec:moreaug} which is hidden in the proof of Proposition 
\ref{pro:seq}.

\begin{figure}
\psfrag{1}[Br][Bl][1.3][0]{$1$}
\psfrag{2}[Br][Bl][1.3][0]{$2$}
\psfrag{3}[Br][Bl][1.3][0]{$3$}
\psfrag{4}[Br][Bl][1.3][0]{$4$}
\psfrag{5}[Br][Bl][1.3][0]{$5$}
\psfrag{6}[Br][Bl][1.3][0]{$6$}
\psfrag{7}[Br][Bl][1.3][0]{$7$}
\psfrag{8}[Br][Bl][1.3][0]{$8$}
\psfrag{9}[Br][Bl][1.3][0]{$9$}
\psfrag{10}[Br][Bl][1.3][0]{$10$}
\psfrag{11}[Br][Bl][1.3][0]{$11$}
\psfrag{12}[Br][Bl][1.3][0]{$12$}
\psfrag{13}[Br][Bl][1.3][0]{$13$}
\psfrag{14}[Br][Bl][1.3][0]{$14$}
\psfrag{15}[Br][Bl][1.3][0]{$15$}
\psfrag{16}[Br][Bl][1.3][0]{$16$}
\psfrag{17}[Br][Bl][1.3][0]{$17$}
\psfrag{18}[Br][Bl][1.3][0]{$18$}
\psfrag{19}[Br][Bl][1.3][0]{$19$}
\psfrag{20}[Br][Bl][1.3][0]{$20$}
\psfrag{21}[Br][Bl][1.3][0]{$21$}
\psfrag{22}[Br][Bl][1.3][0]{$22$}
\psfrag{23}[Br][Bl][1.3][0]{$23$}
\psfrag{24}[Br][Bl][1.3][0]{$24$}
\psfrag{25}[Br][Bl][1.3][0]{$25$}
\psfrag{26}[Br][Bl][1.3][0]{$26$}
\psfrag{1'}[Bl][Bl][1.3][0]{$1$}  
\psfrag{2'}[Bl][Bl][1.3][0]{$2$}  
\psfrag{3'}[Bl][Bl][1.3][0]{$3$}  
\psfrag{4'}[Bl][Bl][1.3][0]{$4$}  
\psfrag{5'}[Bl][Bl][1.3][0]{$5$}
\psfrag{6'}[Bl][Bl][1.3][0]{$6$}
\psfrag{7'}[Bl][Bl][1.3][0]{$7$}
\psfrag{8'}[Bl][Bl][1.3][0]{$8$}
\psfrag{9'}[Bl][Bl][1.3][0]{$9$}
\psfrag{10'}[Bl][Bl][1.3][0]{$10$}
\psfrag{11'}[Bl][Bl][1.3][0]{$11$}
\psfrag{12'}[Bl][Bl][1.3][0]{$12$}
\psfrag{13'}[Bl][Bl][1.3][0]{$13$}
\psfrag{14'}[Bl][Bl][1.3][0]{$14$}   
\psfrag{15'}[Bl][Bl][1.3][0]{$15$}
\psfrag{16'}[Bl][Bl][1.3][0]{$16$}
\psfrag{17'}[Bl][Bl][1.3][0]{$17$}
\psfrag{18'}[Bl][Bl][1.3][0]{$18$}  
\psfrag{19'}[Bl][Bl][1.3][0]{$19$}  
\psfrag{20'}[Bl][Bl][1.3][0]{$20$}
\psfrag{21'}[Bl][Bl][1.3][0]{$21$}
\psfrag{22'}[Bl][Bl][1.3][0]{$22$}
\psfrag{23'}[Bl][Bl][1.3][0]{$23$}
\psfrag{24'}[Bl][Bl][1.3][0]{$24$}
\psfrag{25'}[Bl][Bl][1.3][0]{$25$}
\psfrag{26'}[Bl][Bl][1.3][0]{$26$}
\includegraphics[angle=45,width=\linewidth]{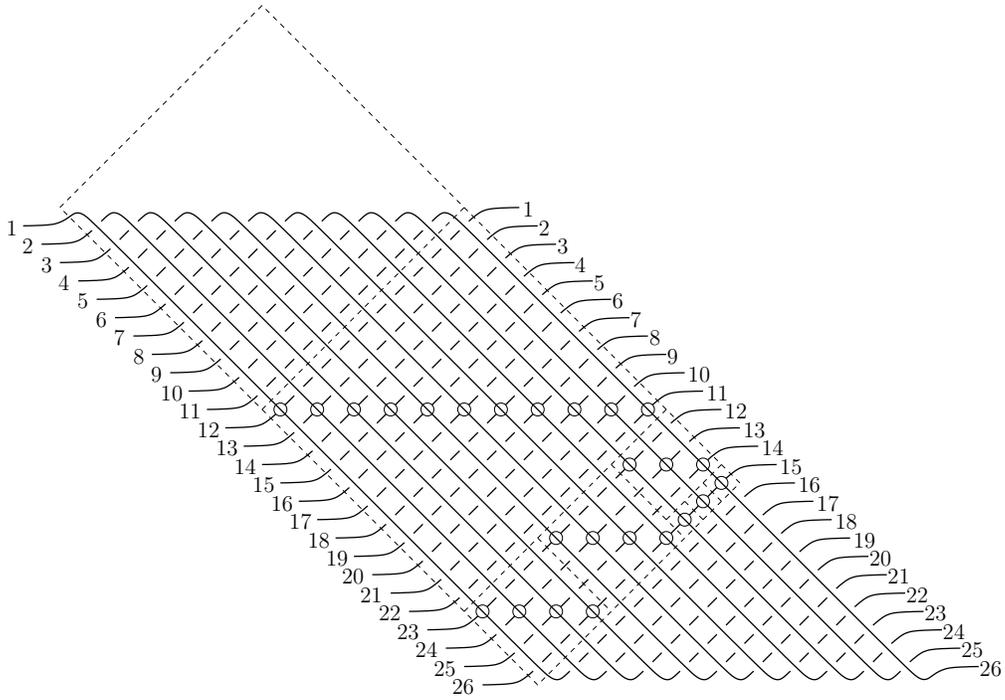}
\caption{The augmentation of a $(p,q)$ torus link 
implements the Euclidean algorithm; on the diagram, 
$p=11$ and $q=26$.}\label{fig:eukl}
\end{figure}

Let us denote the quotients and residues in the Euclidean algorithm 
(with input $p$ and $q$) as follows:
\begin{align*}
p&=k_{-1}q+r_0 & (0&\le r_0<q)\\
q&=k_0r_0+r_1 & (0&\le r_1<r_0)\\
r_0&=k_1r_1+r_2 & (0&\le r_2<r_1)\\
r_1&=k_2r_2+r_3 & (0&\le r_3<r_2)\\
&\vdots&&\\
r_{l-3}&=k_{l-2}r_{l-2}+r_{l-1} & (0&\le r_{l-1}<r_{l-2})\\
r_{l-2}&=k_{l-1}r_{l-1}+r_{l} & (0&\le r_{l}<r_{l-1})\\
r_{l-1}&=k_lr_l+0.&&
\end{align*}
(Of course, $r_l=\gcd\{\,p,q\,\}$.) The points of the augmentation $X$ are 
arranged 
in blocks of the following sizes: $(k_0-1)$ blocks of size $r_0$; $k_1$ 
blocks of size $r_1$; $k_2$ blocks of size $r_2$ and so on until the last 
$k_l$ blocks of size $r_l$. If we draw the diagram of the braid as on Figure 
\ref{fig:eukl}, every block can be viewed as the diagonal of a square, and 
the squares can in turn be seen to be placed inside a $p\times q$ rectangular 
box so that they realize a `graphic implementation' of the 
Euclidean algorithm. 
\end{megj}

\section{A loop of positive links}\label{sec:loop}
\begin{figure}
\resizebox{\linewidth}{506pt}{\includegraphics{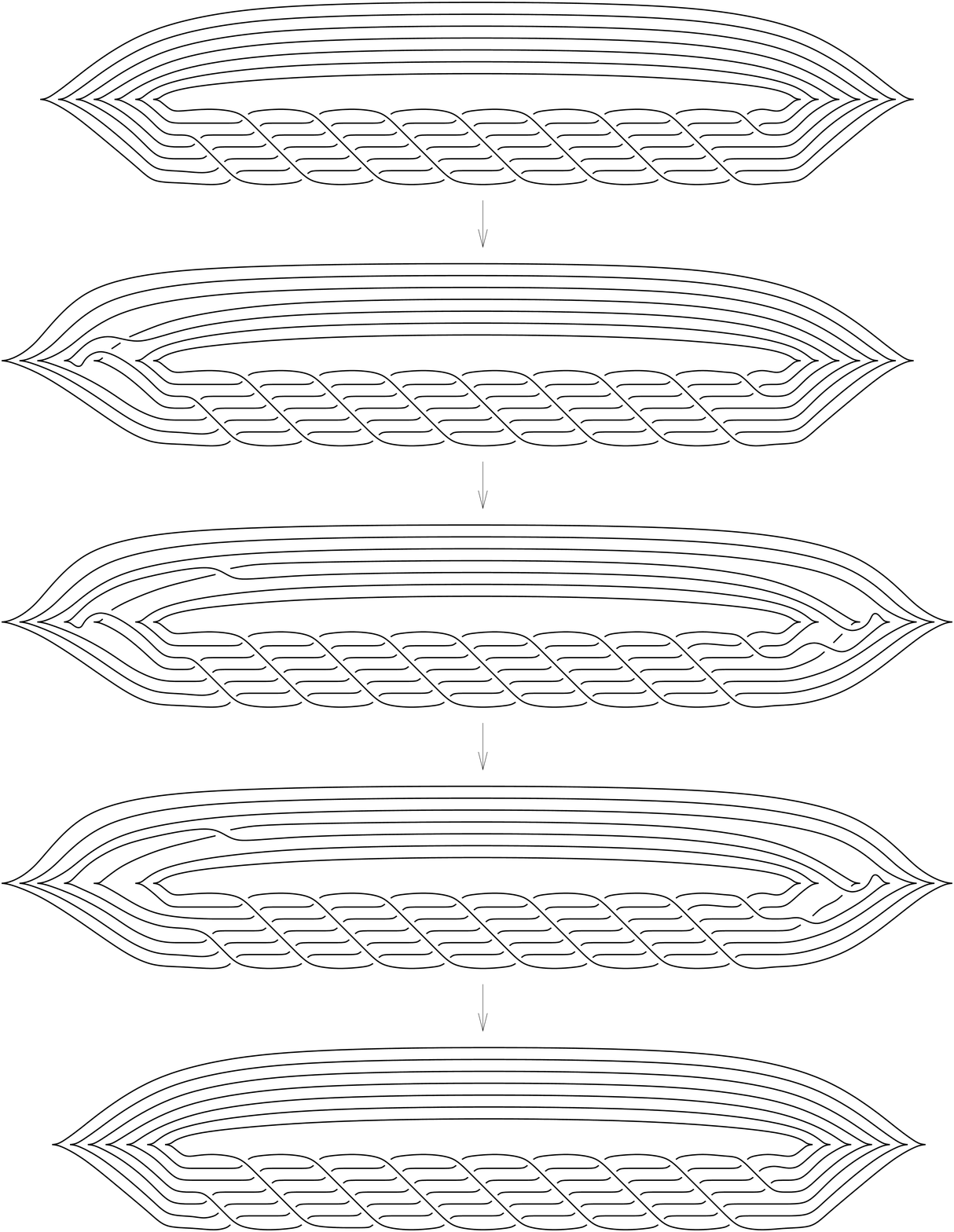}}
\caption{A path that corresponds to conjugating a braid, viewed in the
front projection. The moves are similar to those in Figure 1
of \cite{en}.}\label{fig:97front}
\end{figure}

\begin{figure}
\psfrag{a}[Bl][Bl][4][0]{$a_{\text{temp}}$}
\psfrag{al}[Bl][Bl][4][0]{$a_m$}
\psfrag{all}[Bl][Bl][4][0]{$a_{m+1}$}
\psfrag{b}[Bl][Bl][4][0]{$b_{m,1}$}
\psfrag{c}[Bl][Bl][4][0]{$c_{m,1}$}
\psfrag{e}[Bl][Bl][4][0]{$b_{m+1,1}$}
\psfrag{d}[Bl][Bl][4][0]{$\ddots$}
\psfrag{u}[Bl][Bl][4][0]{$b_{q-1,1}$} 
\psfrag{U}[Bl][Bl][4][0]{$U_{m+2}$}
\resizebox{\linewidth}{492pt}{\includegraphics{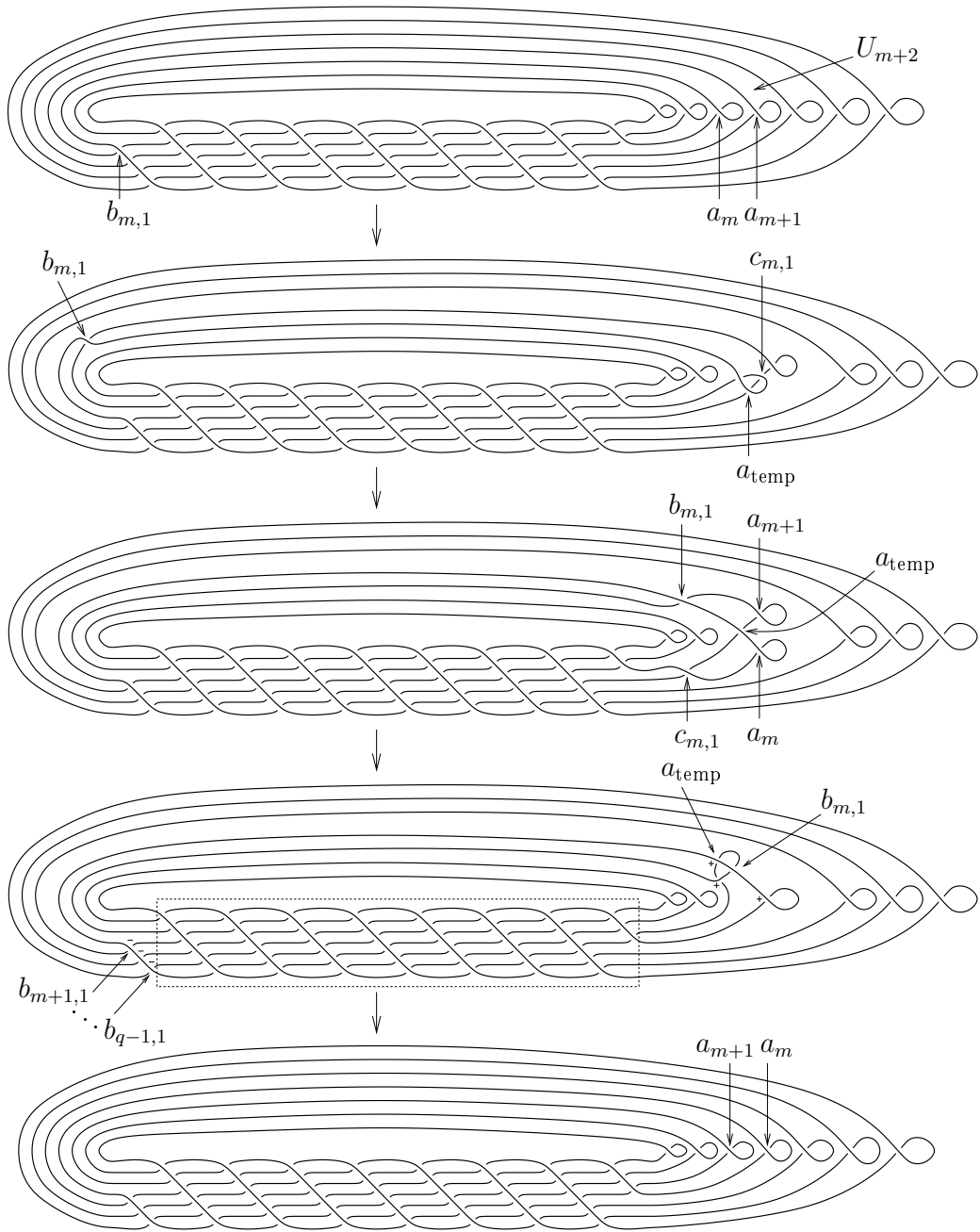}}
\caption{A path 
that corresponds to conjugating a braid, viewed in the Lagrangian 
projection. The individual Reidemeister moves \emph{do not} correspond 
to those in Figure \ref{fig:97front}.}\label{fig:97lag}
\end{figure}

Let $L_\beta$ be the Legendrian closure of the positive braid $\beta$. 
Then there exists a natural closed loop in the connected component 
$\mathscr L_\beta$ of the space of Legendrian links 
that contains $L_\beta$, as follows. Let us write 
$\beta=\lambda_1\ldots\lambda_w$ as a product of the braid group 
generators. On Figure \ref{fig:97front}, we show through an example how 
$L_\beta$ 
can be changed into the Legendrian closure of the conjugate braid that 
results from moving the first factor $\lambda_1$ to the end of the word: 
if $\lambda_1$ is a half-twist of the $m$'th and $(m+1)$'st strands of the 
braid, then one interchanges the $m$'th and $(m+1)$'st strands 
\emph{above} the braid. 
On Figure \ref{fig:97lag}, the same path $\Phi_{\lambda_1}$ is shown, but in 
the Lagrangian projection. Note that the index $1$ crossings $a_m$ and 
$a_{m+1}$ trade places. (The notation used on the diagram for the index 
$0$ crossings is the one that we will introduce below for the special case 
of torus links.) The Lagrangian diagrams of the endpoints are clearly 
obtained by resolution of the corresponding fronts. However, we will not 
prove that the paths themselves agree, too (up to homotopy). 
Instead, we will content ourselves with checking (using Theorem 4.1 
of \cite{en}) that the four Reidemeister moves on Figure \ref{fig:97lag} are 
consistent, and thereafter use the Lagrangian construction as our definition 
of $\Phi_{\lambda_1}$. (There are no such consistency issues with Reidemeister 
moves of fronts, but we need the Lagrangian diagrams to compute holonomies.)

\begin{tetel}\label{thm:consistent}
The sequence of Reidemesiter moves in Figure \ref{fig:97lag}, defining 
$\Phi_{\lambda_1}$, is consistent.
\end{tetel}

\begin{proof}
This is a generalization of Example 4.2
of \cite{en}. Isotope the diagram just like in the proof of 
Lemma \ref{lem:kicsi}. Then by choosing a small enough $\varepsilon$, 
the height $h(a_{m+1})$ will dominate the expression whose positivity 
is needed (by Theorem 4.1
of \cite{en}) in order for the first move (which is a Reidemeister II 
move) to be consistent. To carry out the III$_{\text{b}}$ move that 
follows, we 
need to isotope the second diagram from the top so that what remains from 
the region $U_{m+2}$ after the 
first move has larger area than the vanishing triangle. This can be 
achieved by the 
same trick (this time, moving away the `outer' $q-m-1$ strands). Next, 
the newborn triangle needs to be blown up so that it has larger area than 
the triangle which is due to vanish in the second III$_{\text{b}}$ move 
(which is the third move altogether). 
For this, the same trick with the bulges 
still works: apply it to the $q-m-1$ outer strands and the one 
that crosses itself at $a_m$. Finally, the exact same 
argument guarantees that the fourth move, of type II$^{-1}$, is 
consistent too. 
\end{proof}

Let us use the formulas of section \ref{sec:z2} 
to compute the holonomy $\mu_1$ of $\Phi_{\lambda_1}$ (cf.\ Example 3.4 
of \cite{en}). More 
precisely, we will compute the action of $\mu_1$ on the index $0$ 
crossings which generate the index $0$ contact homology $H_0(L_\beta)$. 
In the first move, two new crossings appear; let us denote the one 
with index $0$ by\footnote{In other words, let it inherit the labels 
of the crossing that is being moved to the other end of the braid. This is 
what we've done in Figure \ref{fig:97lag} too, except that there, a different 
notation is used for index $0$ crossings.} 
$c_{m,\sigma(m+1),1}$ 
and the one with index $1$ by $a_{\text{temp}}$. The old index $0$ 
crossings are not affected by this move 
(i.e., the holonomy maps them identically). 
This is true by Proposition 3.5 
of \cite{en}:
for index reasons, the boundary of any index $0$ crossing is $0$. 
It is easy to see that the following two triangle moves don't affect the old 
index $0$ crossings, either. In the fourth, Reidemeister II$^{-1}$ move 
however the 
crossing $b_{m,\sigma(m+1),1}$ (together with $a_{\text{temp}}$) vanishes, 
and its image in the holonomy 
becomes the polynomial $M'_{m+1,m}=C'_{m+1,m}$. We used primed symbols 
to remind us that these are to be computed with respect to the conjugated 
braid. 
The fact that for each 
admissible disc that turns at the positive 
quadrant at $a_{\text{temp}}$ that faces away from $b_{m,\sigma(m+1),1}$, 
the intermediate sequence of labels 
(of other index $1$ generators through which the boundary of the admissible 
disk passes) doesn't contain 
any number larger than $m$ can be shown just like in the proof of Lemma 
\ref{lem:kicsi}. The fact that the sequence is admissible follows by an 
argument very similar to the proof of Proposition \ref{pro:blank}. 
Finally, the fact that each such admissible sequence does contribute the said 
terms can be established as in the proof of Theorem \ref{thm:relacio}.
We have proven:

\begin{all}\label{pro:mualt}
The holonomy $\mu_1$ of $\Phi_{\lambda_1}$ maps each crossing of the braid 
identically except the first one from the left, which is mapped to the 
polynomial $M'_{m+1,m}=C'_{m+1,m}$ (if the first crossing is between the 
$m$'th and $(m+1)$'st strands). This expression is to be computed as in 
Definition \ref{def:C}, with respect to the conjugate braid 
$\lambda_2\lambda_3\ldots\lambda_w\lambda_1$.
\end{all}

Now, it is clear that the concatenation 
$\Omega_\beta=\Phi_{\lambda_1}\ldots\Phi_{\lambda_w}$ is a closed loop in 
$\mathscr L_\beta$ and by definition, its monodromy is the composition of 
holonomies $\mu_w\circ\ldots\circ\mu_1$, followed by a re-labeling to restore 
the original labels. Namely, each symbol $c$ labeling an index $0$ crossing 
needs to be changed back to $b$; 
in fact, $a_i$ would have to be changed to $a_{\sigma(i)}$, but because 
we only concern ourselves with $H_0(L_\beta)$, this can be ignored.

\section{\boldmath $|\mu_0|$ divides $p+q$}\label{sec:moreBC}
The last two sections of the paper contain the proof of Theorem 
\ref{thm:monorendje}. 
The argument works for any $p,q$ except when $q$ divides $p$ 
(the case of a pure braid), or $p$ divides $q$. The reason why we don't claim 
Theorems \ref{thm:monorendje} and \ref{thm:hurokrendje} for multi-component 
torus links is that in \cite{en}, we only proved Theorems 3.2
and 3.8
(and thus Theorem 1.1)
for knots. (However, the extension of those proofs should only be a matter 
of changing the 
formalism to that of link contact homology.) 

Let us revisit the loop $\Omega_{p,q}$ of Legendrian $(p,q)$ torus 
knots defined in the introduction of \cite{en}. As the braid $\beta$ is 
now composed of $p$ periods, the general loop described in section 
\ref{sec:loop} is the $p$-fold 
concatenation of another, and the latter is easy to identify as 
homotopic\footnote{Because we'll omit the rigorous justification of this 
fact, the reader may 
treat the new description of $\Omega_{p,q}$ as the definition.} 
to $\Omega_{p,q}$. In particular, by Theorem 3.8
of \cite{en}, the `full' monodromy takes the form 
$\mu^p$, where $\mu$ is the monodromy of $\Omega_{p,q}$. From now on, we will 
concentrate on this map $\mu$, and especially on its restriction $\mu_0$ 
to the 
index $0$ part of the contact homology of the standard torus link diagram 
$\gamma$ (shown on Figure \ref{fig:pic}), representing the base point $L$. 

We will adjust our notation to this special situation. The crossings of 
$\beta$ will be indexed with two integers (as opposed to three), namely 
$b_{m,n}$ ($m=1,\ldots,q-1$, $n=1,\ldots,p$) will denote the $m$'th crossing 
counted from the top in the $n$'th period of the braid. 
Note that in the definition of $\mu$, after a full period of the 
braid has been moved from the left end to the right end, a re-labeling 
takes place, too: the second label of each crossing in the other $(p-1)$ 
periods is reduced by $1$, and the labels $c_{m,1}$ in the now last period 
are changed to $b_{m,p}$.

\begin{all}\label{pro:muspec}
The monodromy $\mu$ of the loop $\Omega_{p,q}$ of Legendrian torus links 
acts on the index $0$ generators as follows:
\[\mu(b_{m,n})=\left\{\begin{array}{ll}
b_{m,n-1} & \text{if }2\le n\le p\\
C_{q,m} & \text{if }n=1
\end{array}\right..\]
\end{all}

\begin{proof}
The claim is clear for those crossings not in the first period: by 
Proposition \ref{pro:mualt}, they are only affected, and in the described way, 
by the re-labeling. The rest of the statement will be proven by induction on 
$q-m$. When this value is $1$, i.e.\ $m=q-1$, this is just the statement of 
Proposition \ref{pro:mualt} (the conjugate braid in this case is the 
original $\beta$ again, and the re-labeling changes $C'_{q-1+1,q-1}$ into 
$C_{q,q-1}$). Assume the statement holds for $b_{q-1,1},\ldots, 
b_{m+1,1}$.
Right after the conjugation that removes it from the left end of the braid, 
the image of $b_{m,1}$ is $M''_{m+1,m}$, computed with respect to the braid 
after this conjugation. This can be re-written (by grouping terms with 
respect to the first factor in the product) as
\[M''_{m+1,m}=b_{m+1,1}M'_{m+1,m}+b_{m+2,1}M'_{m+2,m}+\ldots+
b_{q-1,1}M'_{q-1,m}+M'_{q,m},\]
where the terms labeled $M'$ on the right are to be computed in the braid 
indicated 
by the box in Figure \ref{fig:97lag}. Note however that by Definition 
\ref{def:C}, the same expressions are obtained 
if we use the whole braid $\beta$ 
(before re-labeling). If we apply the holonomies of the remaining $q-1-m$ 
conjugations and the re-labeling to this expression, we get (by the inductive
hypothesis)
\[\mu(b_{m,1})=C_{q,m+1}M_{m+1,m}+C_{q,m+2}M_{m+2,m}+\ldots+
C_{q,q-1}M_{q-1,m}+M_{q,m}=C_{q,m}.\]
The last equality is true because the middle expression is exactly what 
results if we group terms in $C_{q,m}$ with respect to the last label 
in the admissible sequence which is more than $m$.
\end{proof}

\begin{all}\label{pro:muBC}
In the contact homology ring $H(L)$, 
we have:
\begin{equation}\label{eq:muB}\mu(B_{i,j})=\left\{\begin{array}{ll}
B_{i-1,j-1}+B_{i-1,q}b_{j-1,p}&\text{if }i,j\ge2\\
B_{i-1,q}&\text{if }i\ge2\text{ and }j=1\\
b_{j-1,p}&\text{if }i=1\text{ and }j\ge2
\end{array}\right.,\end{equation}
\begin{equation}\label{eq:muC}\mu(C_{i,j})=\left\{\begin{array}{ll}
C_{i-1,j-1}&\text{if }j\ge2\\
M_{i-1,q}&\text{if }j=1
\end{array}\right.,\end{equation}
and
\begin{equation}\label{eq:muM}
\mu(M_{i,j})=M_{i-1,j-1},\text{ whenever }i,j\ge2.\end{equation}
\end{all}

We omitted $B_{1,1}$ and $C_{i,i}$ because they (and hence their images) 
are equal to $1$ in the contact homology. Recall also that $M_{i,1}=B_{i,1}$ 
and $M_{1,j}=B_{1,j}$.

\begin{proof}
When $i\ge2$, none of the terms in $B_{i,j}$ contains any of 
$b_{1,1},\ldots,b_{q-1,1}$, so they only need to be re-labeled. This means 
that all crossings that the path through the braid which generated the term 
turned at, are shifted to the 
left by a unit. This operation changes the entry point from the one labeled 
$i$ to the one labeled $i-1$. 
If $j=1$, then the shifted path can be completed 
by the overcrossing strand of the last period of the braid, which shows that 
the re-labeled expression is a summand in $B_{i-1,q}$. Moreover, all such 
summands 
are obtained in this way exactly once. When $j\ge2$, the re-labeling results 
in a summand of $B_{i-1,j-1}$, but not all such are obtained: we miss 
contributions from paths that turn at the crossing (the last one on the 
strand with right endpoint $j-1$) $b_{j-1,p}$. Hence the correction term in 
the top row of \eqref{eq:muB} (note that the paths turning at $b_{j-1,p}$ 
are exactly those that would otherwise have arrived at $q$).

When $i=1$ and $j\ge2$, we have 
\begin{equation*}\begin{split}
\mu(B_{1,j})&=
\mu(b_{1,1}B_{2,j}+b_{2,1}B_{3,j}+\ldots+b_{q-1,1}B_{q,j}+R)\\
&=C_{q,1}(B_{1,j-1}+B_{1,q}b_{j-1,p})+C_{q,2}(B_{2,j-1}+B_{2,q}b_{j-1,p})+
\ldots\\
&+C_{q,q-1}(B_{q-1,j-1}+B_{q-1,q}b_{j-1,p})+\mu(R)\\
&=C_{q,1}B_{1,j-1}+C_{q,2}B_{2,j-1}+\ldots+C_{q,q-1}B_{q-1,j-1}\\
&+(C_{q,1}B_{1,q}+C_{q,2}B_{2,q}+\ldots+C_{q,q-1}B_{q-1,q})b_{j-1,p}\\
&+B_{q,j-1}+B_{q,q}b_{j-1,p}\\
&=C_{q,j-1}+C_{q,j-1}C_{j-1,j-1}+C_{q,q}b_{j-1,p}
=C_{q,j-1}+C_{q,j-1}+b_{j-1,p}\\
&=b_{j-1,p}.
\end{split}\end{equation*}

Here, $R$ is the sum of the contributions to $B_{1,j}$ that don't contain 
crossings of the first period. These terms only have to be re-labeled and that 
can be done just like in the argument above for $B_{i,j}$ when $i\ge2$.
In the sum of sums 
\[C_{q,1}B_{1,j-1}+C_{q,2}B_{2,j-1}+\ldots+C_{q,q-1}B_{q-1,j-1}+B_{q,j-1},\] 
we re-grouped the terms; those with an admissible sequence of labels 
formed $C_{q,j-1}$, and the rest, where $j-1$ was repeated `illegally,' 
formed $C_{q,j-1}C_{j-1,j-1}$.

Note that by the now proven \eqref{eq:muB}, for all $i,j\ge2$, 
$\mu(B_{i,j}+B_{i,1}B_{1,j})=B_{i-1,j-1}$. Therefore, when $i>j\ge2$, 
\begin{equation*}\begin{split}
\mu(C_{i,j})&=
\mu\left(\sum_{\{\,i_1,\ldots,i_c,j\,\}\in D_i}
B_{i,i_1}B_{i_1,i_2}B_{i_2,i_3}\ldots B_{i_{c-1},i_c}B_{i_c,j}\right)\\
&=\mu\left(\sum_{\text{\scriptsize$\begin{gathered}\{\,j_1,\ldots,j_d,j\,\}\in D_i\\
j_1,\ldots,j_d\ge2\end{gathered}$}}
\begin{aligned}[t]
(B_{i,j_1}+B_{i,1}B_{1,j_1})(B_{j_1,j_2}+&B_{j_1,1}B_{1,j_2})\ldots\\
(B_{j_{d-1},j_d}+B_{j_{d-1},1}B_{1,j_d})&(B_{j_d,j}+B_{j_d,1}B_{1,j})
\end{aligned}\right)\\
&=\sum_{\text{\scriptsize$\begin{gathered}\{\,j_1,\ldots,j_d,j\,\}\in D_i\\
j_1,\ldots,j_d\ge2\end{gathered}$}}B_{i-1,j_1-1}B_{j_1-1,j_2-1}\ldots
B_{j_{d-1}-1,j_d-1}B_{j_d-1,j-1}\\
&=C_{i-1,j-1}.
\end{split}\end{equation*}
As a consequence of this and \eqref{eq:muB}, for all $i\ge2$,
\[\mu(C_{i,1})=\mu\left(B_{i,1}+\sum_{j=2}^{i-1}C_{i,j}B_{j,1}\right)
=B_{i-1,q}+\sum_{j=2}^{i-1}C_{i-1,j-1}B_{j-1,q}=M_{i-1,q}.\]
Finally, when $2\le i<j$,
\begin{equation*}\begin{split}
\mu(M_{i,j})&=\mu\left(B_{i,j}+\sum_{k=1}^{i-1} C_{i,k}B_{k,j}\right)\\
&=B_{i-1,j-1}+B_{i-1,q}b_{j-1,p}+M_{i-1,q}b_{j-1,p}\\
&+\sum_{k=2}^{i-1} C_{i-1,k-1}(B_{k-1,j-1}+B_{k-1,q}b_{j-1,p})\\
&=M_{i-1,q}b_{j-1,p}+\left(B_{i-1,j-1}
+\sum_{k=2}^{i-1}C_{i-1,k-1}B_{k-1,j-1}\right)\\
&+\left(B_{i-1,q}+\sum_{k=2}^{i-1}C_{i-1,k-1}B_{k-1,q}\right)b_{j-1,p}\\
&=M_{i-1,q}b_{j-1,p}+M_{i-1,j-1}+M_{i-1,q}b_{j-1,p}\\
&=M_{i-1,j-1},
\end{split}\end{equation*}
and when $i>j$, the argument is very similar to the one we gave for $C_{i,j}$.
\end{proof}

\begin{tetel}\label{thm:szamol}
The order of the (restricted) 
monodromy $\mu_0=\mu\big|_{H_0(L)}$ of the loop $\Omega_{p,q}$ of 
Legendrian $(p,q)$ torus knots divides $p+q$.
\end{tetel}

\begin{proof}
This is a straightforward computation, generalizing the first paragraph 
of the proof of Proposition 3.7 
of \cite{en}. Consider the generator $b_{m,p}$ 
($m=1,\ldots,q-1$). By Proposition \ref{pro:muspec}, 
the first $p$ iterations of $\mu$ act on it as follows:
\[\mu(b_{m,p})=b_{m,p-1},\text{ }\mu^2(b_{m,p})=b_{m,p-2},\dots,
\mu^{p-1}(b_{m,p})=b_{m,1},\text{ }\mu^p(b_{m,p})=C_{q,m}.\]
Then by \eqref{eq:muC} of Proposition \ref{pro:muBC}, 
the next $m$ iterations are as follows:
\[\mu^{p+1}(b_{m,p})=C_{q-1,m-1},\ldots,\mu^{p+m-1}(b_{m,p})=C_{q-m+1,1},
\text{ }\mu^{p+m}(b_{m,p})=M_{q-m,q}.\]
Now by \eqref{eq:muM}, the next $q-m-1$ iterations are
\[\mu^{p+m+1}(b_{m,p})=M_{q-m-1,q-1},\ldots,\mu^{p+q-1}(b_{m,p})=M_{1,m+1}.\]
Finally, because $M_{1,m+1}=B_{1,m+1}$, \eqref{eq:muB} yields 
\[\mu^{p+q}(b_{m,p})=b_{m,p}.\]
Because $b_{m,n}$ is on the orbit of $b_{m,p}$ for all $n=1,\ldots,p$, we see 
that $\mu^{p+q}$ is identical on all of the degree $0$ generators.
\end{proof}

\section{\boldmath $(p+q)$ divides $|\mu_0|$}\label{sec:moreaug}
\begin{tetel}\label{thm:osztja}
The sum $p+q$ divides the order of the monodromy $\mu$. 
\end{tetel}

In the proof of Theorem \ref{thm:szamol}, we described explicitly the 
$(p+q)$-element orbit 
of each index $0$ generator $b_{m,n}$ 
(altogether $q-1$ orbits). Recall that all of those orbits
contain a \emph{$b$--sequence} $b_{m,p},\ldots,b_{m,1}$ of length $p$, a 
\emph{$C$--sequence} $C_{q,m},\ldots,C_{q-m+1,1}$ of length $m$ and an 
\emph{$M$--sequence} $M_{q-m,q},\ldots,M_{1,m+1}$ of length $q-m$.
(Recall also that these expressions are cycles in the chain complex 
$\mathscr A$, but we really mean the homology classes represented by 
them.) It suffices to find a suitable one among these orbits 
and show that it isn't periodic by any period shorter than $p+q$. The way 
we will achieve this is the evaluation of the augmentation 
$\varepsilon=\varepsilon_X$ (see Proposition \ref{pro:aug}) on elements of 
the orbit, and proving that the resulting sequence of $0$'s and $1$'s, 
which we will call the \emph{$0$-$1$--sequence of the orbit}, has 
no such shorter period. In fact, we claim the following, from which Theorem 
\ref{thm:osztja} follows immediately: 

\begin{all}\label{pro:01}
If $q>p$ but $p\nmid q$, the $0$-$1$--sequence $S$ of the orbit of $b_{p,p}$ 
consists of $p$ consecutive $0$'s and $q$ consecutive $1$'s.
If $q<p$ but $q\nmid p$, then the same holds for the orbit of 
$b_{[p\pmod q],p}$.
\end{all}

In the latter case, we will denote the value $1\le[p\pmod q]\le q-1$ by 
$r_0$. For the rest of the section, either this value $r_0$ or $p$, as the 
case may be, should be substituted for $m$ in the formulas for the $b$--, 
$C$--, and $M$--sequences. 
Note that if $q\mid p$, the said orbit doesn't even 
exist (if $q>p$ and $p\mid q$, then its $0$-$1$--sequence consists only 
of $1$'s). 


Recall that $X$ was constructed so that the graph realized by $X$ was 
the augmented graph of the underlying permutation $\sigma$ of the braid. 
Hence, this oriented graph 
$\Gamma_\sigma$ has adjacency matrix $\left[\varepsilon_X(B_{i,j})\right]$ 
and 
therefore it contains all the information we need to evaluate 
the algebra homomorphism $\varepsilon=\varepsilon_X$ on the polynomial 
expressions of the $C$-- and $M$--sequences. We will only need to refer to the 
actual braid in the case of the $b$--sequence. In our situation, 
\[\sigma(i)=[(i-p)\pmod q],\text{ }i=1,\ldots,q.\]
When $q<p$, we could equivalently write $\sigma(i)=[(i-r_0)\pmod q]$.
This explains why our choice of orbit in Proposition \ref{pro:01} is 
reasonable: both in the $C$--sequence and in the $M$--sequence the two 
lower indices are always the endpoints of an edge of $\Gamma_\sigma$, but 
\emph{they are listed in the reverse order}. So when we evaluate 
$\varepsilon$ on these polynomials, what we need to examine 
is whether 
the given edge of the graph can be ``reversed,'' i.e.\ if it is part 
of an oriented 
loop (and how many loops) with an admissible sequence of vertices.

We will re-state and prove Proposition \ref{pro:01} in a more detailed 
version.

\begin{all}\label{pro:seq}
The orbit specified in Proposition \ref{pro:01} contributes 
$0$'s and $1$'s to the sequence $S$ as follows.
\begin{enumerate}[(1)]

\item\label{egy} If $2p\le q$ but $p\nmid q$ (hence in fact $2p<q$), 
then we have 
\begin{multline*}
\overbrace{b_{p,p},\ldots,b_{p,1}}^{p\text{ copies of }1};
\overbrace{C_{q,p},\ldots,C_{q-p+1,1}}^{p\text{ copies of }1};\\
\overbrace{M_{q-p,q},\ldots,M_{q-2p+1,q-p+1}}^{p\text{ copies of }0},
\overbrace{M_{q-2p,q-p},\ldots,M_{1,p+1}}^{q-2p\text{ copies of }1}.
\end{multline*}

\item\label{ket} If $2p>q$ but $q>p$, then the sequence is
\[
\overbrace{b_{p,p},\ldots,b_{p,q-p+1}}^{2p-q\text{ copies of }0},
\overbrace{b_{p,q-p},\ldots,b_{p,1}}^{q-p\text{ copies of }1};
\overbrace{C_{q,p},\ldots,C_{q-p+1,1}}^{p\text{ copies of }1};
\overbrace{M_{q-p,q},\ldots,M_{1,p+1}}^{q-p\text{ copies of }0}.
\]

\item\label{ha} If $q<p$, $q\nmid p$, and $2r_0<q$, then we get 
\begin{multline*}
\overbrace{b_{r_0,p},\ldots,b_{r_0,r_0+1}}^{p-r_0\text{ copies of }0},
\overbrace{b_{r_0,r_0},\ldots,b_{r_0,1}}^{r_0\text{ copies of }1};
\overbrace{C_{q,r_0},\ldots,C_{q-r_0+1,1}}^{r_0\text{ copies of }1};\\
\overbrace{M_{q-r_0,q},\ldots,M_{r_0+1,2r_0+1}}
^{q-2r_0\text{ copies of }1},
\overbrace{M_{r_0,2r_0},\ldots,M_{1,r_0+1}}^{r_0\text{ copies of }0}.
\end{multline*}

\item\label{negy} Finally, if $q<p$, $q\nmid p$, and $2r_0\ge q$, then 
we have 
\begin{multline*}
\overbrace{b_{r_0,p},\ldots,b_{r_0,q-r_0+1}}
^{p-q+r_0\text{ copies of }0},
\overbrace{b_{r_0,q-r_0},\ldots,b_{r_0,1}}^{q-r_0
\text{ copies of }1};\\
\overbrace{C_{q,r_0},\ldots,C_{q-r_0+1,1}}^{r_0\text{ copies of }1};
\overbrace{M_{q-r_0,q},\ldots,M_{1,r_0+1}}^{q-r_0\text{ copies of }0}.
\end{multline*}
\end{enumerate}
\end{all}

\begin{proof}
\textbf{\boldmath $b$--sequence}. In case \eqref{egy}, all numbers 
$1\le j\le p$ are so that 
$j<[(j+p)\pmod q]=j+p$ and $j<[(j-p)\pmod q]$. This means that $j_+=j_-=j$, 
i.e.\ that there is a loop edge attached to $j$ in $\Gamma_\sigma$. It is 
easy to check that in the construction of $X$, the crossing that realizes 
this loop edge is exactly $b_{p,j}$.

In case \eqref{ket}, we similarly find loop edges but only attached to 
the 
numbers $1,\ldots,q-p$. These are realized by the crossings 
$b_{p,1},\ldots,b_{p,q-p}$. For the crossings $b_{p,q-p+1},\ldots,b_{p,p}$, 
we find that their second labels in the system that we used to label 
crossings 
of general braids in section \ref{sec:BC} are also $q-p+1,\ldots,p$. These 
numbers can not be the endpoints of a chord because the numbers preceding 
them 
in the permutation (namely, $1,\ldots,2p-q$) are smaller than them. 
Therefore these crossings are indeed not selected into $X$.

In cases \eqref{ha} and \eqref{negy}, first note that the crossings 
$b_{r_0,q+1},\ldots,b_{r_0,p}$ have third labels greater than $1$ in 
the labeling system of section \ref{sec:BC}, so they never get 
selected into $X$. Neither do $b_{r_0,r_0+1},\ldots,b_{r_0,q}$, 
because their `old' first labels are $r_0+1,\ldots,q$ and these are 
taken to the smaller values $1,\ldots,q-r_0$ by $\sigma$ (i.e., 
they'll never be the startpoint of a chord). After this, the rest of 
the 
$b$--sequence can be sorted out just like in the first two cases.

\textbf{\boldmath $C$--sequence}. Here the claim is that it always contributes 
only $1$'s to $S$. This is because $i>\sigma^{-1}(i)=[(i+p)\pmod q]$ 
implies 
$\varepsilon(C_{i,[(i+p)\pmod q]})=1$. Indeed, since it is preceded in 
the permutation by a smaller number, $i$ can't
be the endpoint of a chord, only of the single edge coming from 
$\sigma^{-1}(i)$.
So there is a unique term in $C_{i,[(i+p)\pmod q]}$ that
contributes $1$ to $\varepsilon(C_{i,[(i+p)\pmod q]})$,
namely the one which, when multiplied by $B_{[(i+p)\pmod q],i}$ on the
right,
produces the term corresponding to the unique loop described in
Lemma \ref{lem:graf-szint}.

\textbf{\boldmath $M$--sequence}. All the four claims in this case follow from 
Lemma \ref{lem:masodik}. Note in particular that if we allowed 
$p\mid q$ in case \eqref{egy}, then $j-p$ would be the second largest 
vertex on $\Gamma_j$ for all $j=p+1,\ldots,q$, and therefore the first 
$p$ elements of the $M$--sequence wouldn't be mapped to $0$ by 
$\varepsilon$. However if $p\nmid q$, then the first number $x$ 
after the 
sequence $j,j-p,j-2p,\ldots$ `wraps around' the circle $\Z_q$ is 
different from $j$. If $x$ is also smaller than $j$ (and this will be 
the case exactly when $j=q,q-1,\ldots,q-[q\pmod p]+1$), then it falls 
between $j$ and $j-p$, so $j-p$ 
is not second largest in $\Gamma_j$. If $x$ is larger than $j$, then 
it is recognized as what we called $\sigma(j_+)$ in section 
\ref{sec:aug}. But then there exists a chord in $\Gamma_\sigma$, 
starting from $j_+$, and ending at the element $j_-$ of 
$\Gamma_j$. This $j_-$ is by construction such that $[(j_-+p)\pmod 
q]>j$. But because $p<q/2$, if we assume that $j>q-p$, 
then this is 
only possible if $j_-+p>j$, i.e.\ if $j_->j-p$. This again means that 
in these cases, $j-p$ is not second largest on $\Gamma_j$. This proves 
the claim about the first part of the $M$--sequence. Finally, if 
$j\le q-p$ (which implies $j_-=j$), then $\Gamma_j$ only visits the 
positive elements of the arithmetic 
progression $j,j-p,j-2p,\ldots$  (i.e., there is no wrapping around) 
and $j-p$ is obviously second largest among these. Cases \eqref{ket}, 
\eqref{ha} and \eqref{negy} can be handled similarly.
\end{proof}

Note that this was a generalization of the middle paragraph in the proof 
of Proposition 3.7
of \cite{en}; there, we used the augmentation $X=\{\,b_3\,\}$.
Thus concludes the proof of Theorem \ref{thm:osztja}, which in turn 
implies 
that the monodromy invariant introduced in \cite{en} is non-trivial.

\end{document}